\documentclass[12pt,twoside]{amsart}
\usepackage{amssymb}
\usepackage{amscd}
\usepackage{xy}

\nonstopmode

\numberwithin{equation}{section}
\font\tengothic=eufm10 scaled\magstep 1
\font\sevengothic=eufm7 scaled\magstep 1
\newfam\gothicfam
         \textfont\gothicfam=\tengothic
         \scriptfont\gothicfam=\sevengothic

\newtheorem{theorem}{Theorem}[section]
\newtheorem{lemma}[theorem]{Lemma}
\newtheorem{proposition}[theorem]{Proposition}
\newtheorem{corollary}[theorem]{Corollary}

\theoremstyle{definition}
\newtheorem{definition}[theorem]{Definition} 
\newtheorem{remark}[theorem]{Remark}
\newtheorem{example}[theorem]{Example}

\newtheorem{notation}[theorem]{Notation}
\newtheorem{question}[theorem]{Question}
\newtheorem{algorithm}[theorem]{Algorithm}

\newcommand\Hom{\operatorname{Hom}}
\newcommand\Ext{\operatorname{Ext}}

\newcommand\rank{\operatorname{rank}}

\newcommand\im{\operatorname{im}}

\newcommand{\proj}[1]
{ \mathchoice
              { {\mathbb P}^{#1} }
              { {\mathbb P}^{#1} }
              { {\mathbb P}^{#1} }
              { {\mathbb P}^{#1} }
            }

\newcommand{\s}{\; | \;}
\newcommand{\mif}{\mbox{if} ~}

\newcommand{\ffi}{\varphi}

\newcommand{\cI}{{\mathcal I}}

\newcommand{\cE}{{\mathcal E}}

\newcommand{\cO}{{\mathcal O}}

\newcommand{\cN}{{\mathcal N}}

\newcommand{\cOP}{{\mathcal O}_{\mathbb{P}^3}}
\newcommand{\acm}{arithmetically Cohen-Macaulay }

\newcommand {\ZZ}{\mathbb{Z}}

\newcommand {\PP}{\mathbb{P}}

\newcommand {\FF}{\mathbb{F}}
\newcommand{\eps}{\varepsilon}

\begin{document}
\title[Tetrahedral curves]{Tetrahedral curves}

\author[J.\ Migliore, U.\ Nagel]{J.\ Migliore${}^*$, U.\ Nagel${}^{+}$}
\address{Department of Mathematics,
           University of Notre Dame,
           Notre Dame, IN 46556,
           USA}
\email{Juan.C.Migliore.1@nd.edu}
\address{Department of Mathematics,
University of Kentucky, 715 Patterson Office Tower,
Lexington, KY 40506-0027, USA}
\email{uwenagel@ms.uky.edu}

\thanks{${}^*$ Part of the work for this paper was done while this
author was sponsored by the National Security Agency under Grant
Number MDA904-03-1-0071. \\
${}^{+}$ This author gratefully acknowledges partial support by 
a Special Faculty Research Fellowship from the University of Kentucky.}


\begin{abstract}
A tetrahedral curve is a space curve whose defining ideal is an
intersection of powers of monomial prime ideals of height two.  It is
supported on a tetrahedral configuration of lines.  Schwartau described
when certain such curves are ACM,  namely he restricted to curves
supported on a certain four of the six lines.  We consider the general
situation.

We first show that starting with an arbitrary tetrahedral curve, there is
a  particular reduction that produces a  smaller tetrahedral curve and
preserves the even liaison class.  We call the curves that are minimal
with respect to this reduction S-minimal curves.  Given a tetrahedral
curve, we describe a simple algorithm (involving only integers) that
computes the S-minimal curve of the corresponding even liaison class; in
the process it determines if the original curve is arithmetically
Cohen-Macaulay or not.  We also describe the minimal free resolution of
an S-minimal curve, using the theory of cellular resolutions.  This
resolution is always linear.  This result allows us to classify the
arithmetically Buchsbaum, non-ACM tetrahedral curves.  More importantly,
it allows us to conclude that an S-minimal curve is minimal in its even
liaison class; that is, the whole even liaison class can be built up
from the S-minimal curve.  Finally, we show that there is a large  set of
S-minimal curves such that each curve corresponds to a smooth point of a
component of the Hilbert scheme and that this component has the expected
dimension.

\end{abstract}


\maketitle

\tableofcontents


    \section{Introduction} \label{intro}

In his Ph.D.\ thesis \cite{Sw}, which was never published, Phil Schwartau
considered certain monomial ideals and the question of whether or not they were
   Cohen-Macaulay.  Specifically, he considered ideals of the form
\[
(X_0,X_1)^{a} \cap (X_1,X_2)^b \cap (X_2,X_3)^c \cap (X_3,X_0)^d
\]
in the ring $k[X_0,X_1,X_2,X_3]$ where $k$ is an algebraically closed field.
These are unmixed ideals defining curves in
$\proj{3}$.  They are supported on a complete intersection (taking $(a,b,c,d) =
(1,1,1,1)$ gives the complete intersection of $X_0X_2$ and $X_1X_3$).
Schwartau  gave a complete classification of the 4-tuples of integers
$(a,b,c,d)$ that define arithmetically Cohen-Macaulay curves (see Theorem
\ref{Sw} and Theorem \ref{thm-corr-sw}).

This paper arose from our desire to make the natural extension of this result
to include the lines defined by $(X_0,X_2)$ and $(X_1,X_3)$.  To simplify the
notation we have changed the variables, so we are considering the ring $R =
k[a,b,c,d]$ and ideals of the form
\[
I = (a,b)^{a_1} \cap (a,c)^{a_2} \cap (a,d)^{a_3} \cap (b,c)^{a_4} \cap
(b,d)^{a_5} \cap (c,d)^{a_6}.
\]
Note that $I$ is an unmixed  monomial ideal.  Since these six lines 
can be viewed as
forming the edges of a tetrahedron, we call such a curve a {\em 
tetrahedral curve}. It is useful to consider the empty set as the {\em trivial curve} defined by $(0, 0, 0, 0, 0, 0)$. 

More than simply considering the question of when such a curve is
arithmetically Cohen-Macaulay, we are interested in getting some idea of the
even liaison class of such a curve.  For background on liaison, see the book
\cite{migliore}; for the most part we will assume the necessary definitions and
basic results on liaison.

Many papers in the literature give classification results of the following
kind: they  consider a particular kind of curve, and ask when two such curves
are linked, or ask for a description of the even liaison class of such a curve
(cf.\ for instance
\cite{bresinsky-huneke},
\cite{geominv}, \cite{dble}, \cite{NNS}).  The work in this paper can be
viewed, in part, as a new contribution to this kind of question.  However,
the main interest comes from the naturalness of these monomial ideals 
themselves,
  from the surprising effectiveness of our reduction procedure, from
the simplicity of the minimal free resolution obtained, and from the 
interplay of
strikingly different techniques and tools used to obtain these results.

The key idea that got this work started was a realization that there is a
simple reduction possible for tetrahedral curves (Proposition \ref{reduce}).
Using the machinery of basic double links, this reduction accomplishes two
amazing things: the new curve is in the same even liaison class as the original
curve, and the new curve is again a tetrahedral curve, but smaller!

Using this reduction, one of two things happens:  either the process of
reduction continues until the curve vanishes, i.e.\ we get the trivial curve, or the process stops at a curve
that cannot be further reduced in this way.  The curves that ultimately reduce
to the trivial curve are clearly arithmetically Cohen-Macaulay, thanks to
liaison theory.  Curves that cannot be further reduced in this way are very
important, and we call them {\em S-minimal curves}.  We give a numerical
criterion for S-minimality (Corollary
\ref{s-min} and Lemma \ref{lem-max-fac}).  This leads to a simple algorithm
(easily done by hand) for computing an S-minimal curve starting with any
tetrahedral curve. A MAPLE implementation of it is documented in the appendix. 

The first main result of this paper is Theorem \ref{thm-res}, where the minimal
free resolution of an S-minimal tetrahedral curve is computed.  We show, in
particular, that this resolution is linear.  This proof starts by finding the
minimal generators of the ideal, and then it translates these generators into a
cell complex and uses the theory of cellular resolutions as developed in
\cite{BSt} to find the rest of the resolution.

An immediate consequence of this result is that an S-minimal curve is
arithmetically Cohen-Macaulay if and only if it is trivial.  This allows one
to easily determine, for any 6-tuple of integers $(a_1,a_2,a_3,a_4,a_5,a_6)$,
whether it corresponds to an arithmetically  Cohen-Macaulay curve or not.  
Another easy consequence is that a tetrahedral
curve is arithmetically Buchsbaum if and only if its Hartshorne-Rao module has
diameter 1.  (It was already known by Schwartau that certain arithmetically
Buchsbaum tetrahedral curves existed.)

A much deeper consequence of Theorem \ref{thm-res} is the second main result of
the paper (Theorem \ref{prop-min}), which says that a tetrahedral curve $C$ is
S-minimal if and only if it is minimal in its even liaison class.

We give two applications of this result.  The first, Theorem \ref{thm-corr-sw},
is a new proof of Schwartau's theorem classifying the arithmetically
Cohen-Macaulay curves supported on the complete intersection by giving the
precise set of 4-tuples.  The second application, Corollary \ref{cor-abm}, is a
classification of the 6-tuples that define minimal tetrahedral curves that are
arithmetically Buchsbaum.  Both of these applications primarily use the earlier
results described above, but need the minimality in the even liaison class to
complete the proof.

Our third main line of investigation concerns the question of 
unobstructedness and
the search for nice components of the Hilbert scheme.  (A component 
is said to be
``nice'' if it has the expected dimension $4 \cdot \deg C$.)  A great 
deal of work
has been done concerning nice components of the Hilbert scheme -- we refer the
reader to \cite{Dolcetti} and to \cite{kleppe}, both for important 
results along
these lines, and also for numerous references to other work.
Using the method of Dolcetti in \cite{Dolcetti}, we first prove that 
any curve in
$\mathbb P^3$ with linear resolution and Hartshorne-Rao module of 
diameter $\leq 2$
is unobstructed, and its Hilbert scheme has the expected dimension $4 
\cdot \deg C$
(Proposition~\ref{lem-dolcet}).  We then characterize the 6-tuples 
corresponding to
minimal tetrahedral curves with Hartshorne-Rao modules of diameter $\leq 2$
(Corollary~\ref{cor-abm} and Lemma~\ref{lem-diam-two}), which then are
unobstructed.  Turning to a different kind of minimal tetrahedral 
curve, we show
that a tetrahedral curve defined by $(a_1,0,0,0,0,a_6)$ is also unobstructed.
Since we have shown the unobstructedness of a large class of minimal 
tetrahedral
curves, and have shown experimentally on the computer that others are also
unobstructed, we end with the question of whether in fact all minimal 
tetrahedral
curves are unobstructed.

We end the paper with a number of questions that arise naturally from our work.


\section{Background}

Let $R = k[a,b,c,d]$, where $k$ is a field. We abbreviate by ACM the
term ``arithmetically Cohen-Macaulay.''  Recall that a projective subscheme $V$
is ACM if and only if the deficiency modules all vanish:
\[
H^i(\proj{n}, {\mathcal I}_V (t)) = 0 \ \ \hbox{ for all $t \in {\mathbb Z}$
and all $1 \leq i \leq \dim V$}.
\]

We recall that the notion of {\em basic double linkage}
was introduced by Lazarsfeld and Rao \cite{LR} as a way of adjoining a plane
curve, or more generally a complete intersection, to a given curve $C$ in such
a way as to preserve the even liaison class of $C$.  More precisely,

\begin{definition}
Let $I \subset R$ be the saturated ideal of an unmixed curve $C$ in
$\proj{3}$.  Let $F \in I$ and $G \in R$ be homogeneous polynomials such that
$(F,G)$ is a regular sequence.  Then the ideal $G\cdot I + (F)$ is the
saturated ideal of a curve $Y$ which is linked to $C$ in two steps (i.e.\ is
{\em bilinked} to $C$).  $Y$ is said to be a {\em basic double link} of $C$.
As sets, $Y$ is the union of $C$ and the complete intersection defined by
$(F,G)$.  The degree of $Y$ is $\deg C + (\deg F)(\deg G)$.
\end{definition}

The notion of basic double linkage has generalizations in several different
directions (higher dimension, higher projective spaces, Gorenstein liaison).
We refer to \cite{migliore} for the details.  If one has a way of recognizing
an ideal as being of the above form, then one can replace that ideal by the
simpler ideal $I$, knowing that the liaison class is preserved.  For
instance, if the original ideal is Cohen-Macaulay then so is the new, simpler
ideal.  This is the approach we take below.

Recall that for an unmixed ideal $I$, the {\em $n$-th symbolic power} of $I$,
denoted $I^{(n)}$, is the saturation of the top dimensional part of the ideal
$I^n$.  Recall also that if $I$ is a complete intersection then $I^{(n)} =
I^n$.

\begin{definition}\label{sdef}
Let $T$ be any union of six lines $\ell_1,\dots,\ell_6$ forming the edges of a
``tetrahedron.'' A {\em tetrahedral curve} is the non-reduced scheme $C$
supported on $T$ and defined by  the saturated ideal
$I_C = I_{\ell_1}^{a_1} \cap \dots \cap I_{\ell_6}^{a_6}$, where $a_i \geq 0$
for all $i$.  The number $a_i$ is called the {\em weight}
of the line $l_i$.
\end{definition}

\begin{remark}
If $C$ is a tetrahedral curve then we can perform a change of variables and
obtain a monomial ideal which has a primary decomposition of the form
\[
I = (a,b)^{a_1} \cap (a,c)^{a_2} \cap (a,d)^{a_3} \cap (b,c)^{a_4} \cap
(b,d)^{a_5} \cap (c,d)^{a_6}.
\]
For the remainder of this paper we will assume that a tetrahedral curve is a
monomial ideal of this form.
\end{remark}

In his thesis \cite{Sw}, Schwartau considered the curves with ideals of the
form
\[
(X_0,X_1)^{a} \cap (X_1,X_2)^b \cap (X_2,X_3)^c \cap (X_3,X_0)^d
\]
(in his notation).  This is equivalent to taking $a_2 = a_5 =0$ in Definition
\ref{sdef}.
Schwartau was primarily interested in the question of when these curves are
arithmetically Cohen-Macaulay.  His main result on this problem is the
following, which we now translate to our language.

\begin{theorem}[\cite{Sw}] \label{Sw}
The ideal
\[
(a,b)^{a_1} \cap (a,c)^{0} \cap (a,d)^{a_3} \cap (b,c)^{a_4} \cap
(b,d)^{0} \cap (c,d)^{a_6}
\]
defines an \acm curve in $\proj{3}$ if and only if
\begin{tabbing}
\underline{Case 1}.xxx \=blaxxxxxxxxxxxxxxxxxxxxxh\ \ \ \  \=blah\kill
\underline{Case 1}. \> $a_1,a_3,a_4,a_6 >0$: \> $a_1+a_6 = a_3+a_4+
\epsilon, \hbox{ for } \epsilon = -1,0,1$. \\
\underline{Case 2}. \> $a_1,a_4,a_6 >0, \ a_3 = 0$: \> $a_1+a_6 \leq a_4+1$.
\\
\underline{Case 3A}. \> $a_1,a_4>0, \ a_3 = a_6 = 0$: \> always \\
\underline{Case 3B}. \> $a_1,a_6 >0, \ a_3 = a_4 = 0$: \> never \\
\underline{Case 4}. \> $a_1 >0, \ a_3 = a_4 = a_6 = 0$: \> always \\
\end{tabbing}
\end{theorem}

\begin{remark} \label{what phil meant}
Clearly Schwartau intended some kind of reduction of cases, since for instance
the case $a_2 >0$, $a_1 = a_3 = a_4 = 0$ is not included in his theorem.  We
have given an ``invariant'' version in Theorem \ref{thm-corr-sw} below, which
we think reflects Schwartau's intention.  We also give a new proof.
\end{remark}


\section{Reduction and S-minimality}

The key to our approach to this problem is the following reduction method.

\begin{proposition} \label{reduce}
Let $I = (a,b)^{a_1} \cap (a,c)^{a_2} \cap (a,d)^{a_3} \cap
(b,c)^{a_4} \cap (b,d)^{a_5} \cap (c,d)^{a_6}$ where not all exponents $a_i$
are zero.  Consider the following systems of inequalities:
\[
\begin{array}{ccc}
\begin{array}{rrcl}
(A): & a_1+a_2 & \geq & a_4,\\
& a_1+a_3 & \geq & a_5, \\
& a_2+a_3 & \geq & a_6
\end{array}
&&
\begin{array}{rrcl}
(B): & a_1+a_4 & \geq & a_2,\\
& a_1+a_5 & \geq & a_3, \\
& a_4+a_5 & \geq & a_6
\end{array}
\\
\begin{array}{rrcl}
(C): & a_2+a_4 & \geq & a_1,\\
& a_2+a_6 & \geq & a_3, \\
& a_4+a_6 & \geq & a_5
\end{array}
& \hbox{\hskip 1cm} &
\begin{array}{rrcl}
(D): & a_3+a_5 & \geq & a_1,\\
& a_3+a_6 & \geq & a_2, \\
& a_5+a_6 & \geq & a_4.
\end{array}
\end{array}
\]
For $1 \leq i \leq 6$ let $a_i' = \max \{ 0, a_i -1 \}$.  Then we have
\begin{itemize}
\item[(i)] $(A) \Leftrightarrow$ $I$ is a basic double link of
\[
(a,b)^{a_1'} \cap (a,c)^{a_2'} \cap
(a,d)^{a_3'} \cap (b,c)^{a_4} \cap (b,d)^{a_5} \cap (c,d)^{a_6}
\]
using $F = b^{a_1}c^{a_2}d^{a_3}$ and $G = a$.
\item[(ii)] $(B) \Leftrightarrow$ $I$ is a basic double link of
\[
(a,b)^{a_1'} \cap (a,c)^{a_2} \cap (a,d)^{a_3} \cap (b,c)^{a_4'} \cap
(b,d)^{a_5'} \cap (c,d)^{a_6}.
\]
using $F = a^{a_1}c^{a_4}d^{a_5}$ and $G = b$.
\item[(iii)] $(C) \Leftrightarrow$ $I$ is a basic double link of
\[
(a,b)^{a_1} \cap (a,c)^{a_2'} \cap (a,d)^{a_3} \cap (b,c)^{a_4'} \cap
(b,d)^{a_5} \cap (c,d)^{a_6'}.
\]
using $F = a^{a_2}b^{a_4}d^{a_6}$ and $G = c$.
\item[(iv)] $(D) \Leftrightarrow$ $I$ is a basic double link of
\[
(a,b)^{a_1} \cap (a,c)^{a_2} \cap (a,d)^{a_3'} \cap (b,c)^{a_4} \cap
(b,d)^{a_5'} \cap (c,d)^{a_6'}.
\]
using $F = a^{a_3}b^{a_5}c^{a_6}$ and $G = d$.
\end{itemize}

\end{proposition}

\begin{proof}
We will
prove (i); of course the others are proved similarly. First note that
\[
a \cdot (a,b)^{n-1} + (b^n) = (a,b)^n
\]
for $n \geq 1$.
Now consider the monomial $F = b^{a_1}c^{a_2}d^{a_3}$.  Notice that
\begin{equation} \label{first}
F \in (a,b)^{a_1} \cap (a,c)^{a_2} \cap (a,d)^{a_3} \subset
(a,b)^{a_1'} \cap (a,c)^{a_2'} \cap (a,d)^{a_3'}
\end{equation}
even when one or more of the $a_i =0$.  The three inequalities are
equivalent to
\begin{equation} \label{second}
F \in (b,c)^{a_4} \cap (b,d)^{a_5} \cap (c,d)^{a_6}.
\end{equation}
Hence
\[
F \in (a,b)^{a_1'} \cap (a,c)^{a_2'} \cap
(a,d)^{a_3'} \cap (b,c)^{a_4} \cap (b,d)^{a_5} \cap (c,d)^{a_6}.
\]
     Notice that $(a,F)$ is a
regular sequence.  Hence we can construct a basic double link of the form
\[
J = a \cdot
\left [
(a,b)^{a_1'} \cap (a,c)^{a_2'} \cap
(a,d)^{a_3'} \cap (b,c)^{a_4} \cap (b,d)^{a_5} \cap (c,d)^{a_6}
\right ]
    + (b^{a_1}c^{a_2}d^{a_3}).
\]
We have to show that
\[
J = I =
(a,b)^{a_1} \cap (a,c)^{a_2} \cap
(a,d)^{a_3} \cap (b,c)^{a_4} \cap (b,d)^{a_5} \cap (c,d)^{a_6}.
\]
    The ideal
$J$ is a saturated, unmixed ideal, by the theory of basic double linkage (cf.\
\cite{migliore}), so it is enough to show that $J \subset I$ and that they
define schemes of the same degree.

For the first one, the fact that
\[
\begin{array}{c}
a \cdot \left [ (a,b)^{a_1'} \cap (a,c)^{a_2'} \cap
(a,d)^{a_3'} \cap (b,c)^{a_4} \cap (b,d)^{a_5} \cap (c,d)^{a_6} \right ] \\
\\
\subseteq (a,b)^{a_1} \cap (a,c)^{a_2} \cap
(a,d)^{a_3} \cap (b,c)^{a_4} \cap (b,d)^{a_5} \cap (c,d)^{a_6}
\end{array}
\]
is clear, while the fact that $b^{a_1}c^{a_2}d^{a_3} \in I$ comes from
(\ref{first}) and (\ref{second}) above.

For the degree computation, recall that $\binom{n}{2} + n = \binom{n+1}{2}$.
Then, again using the theory of basic double linkage, the degree of $J$ is
\[
\begin{array}{rcl}
\deg J & = &\displaystyle \binom{a_1}{2} + \binom{a_2}{2} + \binom{a_3}{2} +
\binom{a_4+1}{2} + \binom{a_5+1}{2} + \binom{a_6+1}{2} + (a_1+a_2+a_3) \\
&= &\displaystyle \binom{a_1+1}{2} + \binom{a_2+1}{2} + \binom{a_3+1}{2} +
\binom{a_4+1}{2} + \binom{a_5+1}{2} + \binom{a_6+1}{2} \\
& = & \deg I
\end{array}
\]
as desired.

The converse is immediate: the fact that $I$ is a basic double link as stated
implies the inequalities of (A) by again using(\ref{first}) and
(\ref{second}).
\end{proof}

\begin{notation}
For the rest of this paper we will abbreviate the monomial
ideal
\[
(a,b)^{a_1} \cap (a,c)^{a_2} \cap (a,d)^{a_3} \cap (b,c)^{a_4} \cap
(b,d)^{a_5} \cap (c,d)^{a_6}
\]
(or the corresponding curve) by the 6-tuple $(a_1,a_2,a_3,a_4, a_5,a_6)$.
\end{notation}

\begin{remark}
Basic double linkage is a special case of Schwartau's ``liaison addition.''
(See \cite{Sw} for the original and \cite{GM4} for a generalization.)
Without entering into details, we remark that many interesting tetrahedral
curves arise as liaison additions.  For instance, let $I$ be the ideal of the
six lines, i.e.\ the tetrahedral curve $(1,1,1,1,1,1)$.  Let $F$ be the
polynomial $abcd$ giving the four faces of the ``tetrahedron'' defined by
$I$.  Note that $F$ is double along each of the six lines.  Let $G$ be a
generally chosen cubic in $I$.  Then $(F,G)$ self-links $I$ (as can be see
geometrically, using Bezout's theorem) and one can check that in fact the
$d$-th symbolic power of $I$, $I^{(d)}$, can be expressed as
\[
I^{(d)} = G^{d-1} \cdot I + F \cdot I^{(d-2)}.
\]
Hence we see directly the curves $(d, d, d, d, d, d)$ are ACM. Note however, that the $d$-th power $I^d$ is not saturated, and has embedded points. 
   
Schwartau also obtains many Buchsbaum curves using liaison addition.  We will
return to Buchsbaum curves below.
\end{remark}

\begin{definition}
We say that a tetrahedral curve is {\em S-minimal} if there is no reduction of
the type described in Proposition \ref{reduce}.
\end{definition}

\begin{corollary} \label{s-min}
Consider a tetrahedral curve $C = (a_1,a_2,a_3,a_4,a_5,a_6)$ where
not all $a_i$
are 0.  Assume without loss of generality that $a_6 = \max \{ a_1,\dots,a_6
\}$.    Then $C$ is S-minimal if and only if
\[
\begin{array}{rcl}
a_1 & > & \max \{ a_3 + a_5,a_2+a_4\} \hbox{ and} \\
a_6 & > & \max \{ a_4+a_5,a_2+a_3 \}
\end{array}
\]
\end{corollary}

\begin{proof}
It is immediate to check that if the stated conditions are satisfied then
each of (A), (B), (C) and (D) in Proposition \ref{reduce} has at least one
inequality that is not satisfied, hence $C$ cannot be reduced via Proposition
\ref{reduce}.

Conversely, suppose that $C$ is S-minimal.    The second and third
inequalities of (C) and (D) in Proposition \ref{reduce}are forced to be true by
the assumption that $a_6 =  \max \{ a_1,\dots,a_6 \}$, so S-minimal
implies that
$a_1 > \max \{ a_3 + a_5,  a_2+a_4\}$.  In particular, $a_1 >a_2$, $a_1 > a_3$,
$a_1 > a_4$ and $a_1 >  a_5$, so the first two inequalities of (A) and of (B)
must be true.  Hence again the assumption that $C$ is S-minimal forces $a_6 >
\max \{ a_4+a_5,a_2+a_3 \}$.
\end{proof}

\begin{example}
Even with the assumption that $a_6 =  \max \{ a_1,\dots,a_6 \}$,
the first inequality of Corollary \ref{s-min} alone does not imply
S-minimality, as shown by the example $C = (4,2,2,1,1,4)$, which can be
reduced using (A) to the curve $C' = (3,1,1,1,1,4)$.  This new curve $C'$ is
    S-minimal.  It can be checked that $C'$ (and hence $C$) is not
arithmetically Cohen-Macaulay (see Corollary \ref{cor-acm} as well).
\end{example}


Our next goal is to provide an algorithm that produces an $S$-minimal curve to
a given tetrahedral curve $C$.  To this end we will extend the
definition of weights. We will use the tetrahedron $T = T(C)$ whose edges are
part of the (potentially) supporting lines of the curve $C$.

\begin{definition} \label{def-weights}
The {\em weight of a facet} of $T = T(C)$ is the sum of the weights of the
edges forming its boundary. Similarly, the weight of a pair of skew lines being
a subset of the six given lines is the sum of the weights of the lines.
\end{definition}

Using these concepts we can make Corollary \ref{s-min} more precise. Note that
each of the reductions (A) - (D) in Proposition \ref{reduce} reduces the
weights of the edges of a facet of the tetrahedron. Thus, we will say that we
{\em reduce a facet} if we apply the corresponding reduction.

The next result says in particular that a tetrahedral curve is not minimal if
and only if we can reduce a facet of maximal weight.

\begin{lemma} \label{lem-max-fac}
Let $C$ be a tetrahedral curve $C = (a_1,a_2,a_3,a_4,a_5,a_6)$ where we assume
without loss of generality that $a_6 = \max \{ a_1,\dots,a_6 \} > 0$. Let $w$
be the maximal weight of a facet of the tetrahedron $T = T(C)$. Then the
following conditions are equivalent:
\begin{itemize}
\item[(a)] $C$ is not $S$-minimal.
\item[(b)] $a_1 + a_6 \leq w$.
\item[(c)] One can reduce any of the facets of $T$ having maximal weight $w$.
\end{itemize}
\end{lemma}

\begin{proof}
We begin by showing the equivalence of (a) and (b). According to Corollary
\ref{s-min} we know that $C$ is not $S$-minimal if and only if
\[
\begin{array}{rcl}
a_1 & \leq & \max \{ a_3 + a_5,a_2+a_4\} \quad \hbox{ or} \\
a_6 & \leq  & \max \{ a_4+a_5,a_2+a_3 \}.
\end{array}
\]
But this is equivalent to the condition
$$
a_1 + a_6 \leq \max \{a_1 + a_4+a_5,a_1 + a_2+a_3,  a_3 + a_5 + a_6, a_2+a_4 +
a_6\} = w
$$
as claimed.

Since (c) trivially implies (a), it remains to show that (c) is a consequence
of (b). To this end we distinguish four cases.

{\it Case 1}: Let $w = a_3 + a_5 + a_6$. Then (b) provides $a_3 + a_5 \geq a_1$
showing that we can apply reduction (D).

{\it Case 2}: Let $w = a_2 + a_4 + a_6$. Then we conclude as above that we can
use reduction (C).

{\it Case 3}: Let $w = a_1 + a_2 + a_3$. Then (b) implies $a_2 + a_3 \geq a_6$.
Moreover, we have by assumption $w \geq a_2 + a_4 + a_6$. It implies
$$
a_1 + a_3 \geq a_4 + a_6 \geq a_6 \geq a_5.
$$
Similarly $w \geq a_3 + a_5 + a_6$ provides $a_1 + a_2 \geq a_4$. Thus, we have
shown that we can apply reduction (A).

{\it Case 4}: Let $w = a_1 + a_4 + a_5$. Then we see as in Case 3 that we can
use reduction (B).
\end{proof}

The last result leads to an algorithm for producing $S$-minimal curves that
works with weights only.

\begin{algorithm}[for computing $S$-minimal curves] \label{alg} \mbox{ } \\
Input: $(a_1,\ldots,a_6) \in \ZZ^6$ where all $a_i \geq 0$, the weight vector
of a tetrahedral curve $C$. \\ Output: $(a_1,\ldots,a_6)$, the weight vector of
an $S$-minimal curve obtained by reducing $C$.
\begin{itemize}
\item[1.] compute $i$ such that $a_i = \max \{a_1,\ldots,a_6 \}$ \\[1ex]
\item[2.] {\bf if $a_i = 0$ then return} $(0,\ldots,0)$ \\[1ex]
\item[3.] determine a facet $F$ of maximal weight $w$\\[1ex]
\item[4.] {\bf if $a_i + a_{7-i} > w$ then return} $(a_1,\ldots,a_6)$. \\[1ex]
\item[5.] apply the reduction corresponding to the facet $F$ and {\bf go} to
Step 1 \\[1ex]
\end{itemize}
\end{algorithm}

\begin{remark} \label{rem-alg} (i)
Strictly speaking the scheme above is not an algorithm since it leaves choices
in Steps 1 and 3 in case there is more than one line or facet of maximal
weight. But this can be fixed easily, e.g., by using the 
lexicographic order.  See
the appendix.

(ii) Correctness of Algorithm \ref{alg} follows by Lemma \ref{lem-max-fac}.
Note that the lines with index $i$ and $7-i$ do not intersect.
\end{remark}


\section{The Minimal Free Resolution of an S-minimal Curve}

Notice that the $6$-tuple $(0,\ldots,0)$ corresponds to the ring $R$. It turns
out that it is useful to consider $R$ formally as a curve as we did in
Definition \ref{sdef}. We give it a name.

\begin{definition} \label{def-triv}
The {\it trivial} tetrahedral curve is defined by $(0,\ldots,0)$.
\end{definition}

The following is one of the two main results of this paper.

\begin{theorem} \label{thm-res}
Every non-trivial {\em S-minimal} tetrahedral curve has a linear minimal free
resolution.

More precisely, if the curve $C$ is defined by $(a_1,a_2,a_3,a_4,a_5,a_6)$ and
$a_6 = \max \{a_i\} > 0$ then its minimal free resolution has the form
\begin{equation*}
   0 \to R^{\beta_3}(-a_1-a_6-2)  \to R^{\beta_2}(-a_1-a_6-1)  \to
R^{\beta_1}(-a_1-a_6) \to I_C \to 0
\end{equation*}
where
\begin{eqnarray*}
\beta_1 & = & (a_1 + 1) (a_6 + 1) - \sum_{i=2}^5 \frac{a_i (a_i + 1)}{2} \\
\beta_2 & = & 2 a_1 a_6 + a_1 + a_6  - \sum_{i=2}^5 a_i (a_i + 1) \\
\beta_3 & = & a_1 a_6  - \sum_{i=2}^5 \frac{a_i (a_i + 1)}{2}. \\
\end{eqnarray*}
\end{theorem}

\begin{proof}
Let $C$ be an S-minimal tetrahedral curve, and assume that $C$ is the curve
defined by the tuple $(a_1,a_2,a_3,a_4,a_5,a_6)$.  Without loss of
generality assume that $a_6$ is the largest entry.  Note that $(a,b)^{a_1}$
and $(c,d)^{a_6}$ are disjoint ACM curves, so their intersection is given by
their product (cf.\ \cite{serre}, Corollaire, p.\ 143 or \cite{MM}, 
Theorem 1.2).
In particular, the curve
$(a_1,0,0,0,0,a_6)$ is given by the entries of the following $(a_1 + 
1) \times (a_6
+1)$ matrix:

\begin{equation} \label{bigmatrix}
\left [
\begin{array}{cccccccccccccccccccccc}
a^{a_1}c^{a_6} & a^{a_1} c^{a_6-1}d  & \dots  & a^{a_1}cd^{a_6-1} &
a^{a_1}d^{a_6} \\
a^{a_1-1}bc^{a_6} & a^{a_1-1}b c^{a_6-1}d & \dots  &
a^{a_1-1}bcd^{a_6-1} & a^{a_1-1}bd^{a_6} \\
\vdots & \vdots & & \vdots & \vdots \\
ab^{a_1-1}c^{a_6} & ab^{a_1-1} c^{a_6-1}d & \dots & ab^{a_1-1}cd^{a_6-1} &
ab^{a_1-1}d^{a_6} \\
b^{a_1}c^{a_6} & b^{a_1} c^{a_6-1}d & \dots & b^{a_1}cd^{a_6-1} &
b^{a_1}d^{a_6} \\
\end{array}
\right ]
\end{equation}
By Corollary \ref{s-min}, we have the inequalities
\begin{equation} \label{cond for min}
\begin{array}{rcl}
a_1 & > & a_3 + a_5 \\
a_1 & > & a_2 + a_4 \\
a_6 & > & a_4 + a_5 \\
a_6 & > & a_2 + a_3
\end{array}
\end{equation}
In particular, $a_1 > 0$.

Now we want to describe the ideal $(a_1,a_2,a_3,a_4,a_5,a_6)$.  Note that as
ideals,
\[
(a_1,a_2,a_3,a_4,a_5,a_6) \subset (a_1,0,0,0,0,a_6).
\]
Consider the monomials obtained by deleting the following from the matrix
(\ref{bigmatrix}):
$a_2$ diagonals from the ``Southeast corner,'' $a_3$ diagonals from the
``Southwest corner,'' $a_4$ diagonals from the ``Northeast corner,'' and $a_5$
diagonals from the ``Northwest corner."  We obtain the following shape:

\newsavebox{\build}
\savebox{\build}(200,85)[tl]
{
\begin{picture}(200,125)
\put (0,-7){\line (0,1){100}}

\put (0,-7){\line (1,0){140}}
\put (0,93){\line (1,0){140}}
\put (140,-7){\line (0,1){100}}

\put (0,66){\line (1,1){27}}
\put (0,20){\line (1,-1){27}}
\put (140,66){\line (-1,1){27}}
\put (140,20){\line (-1,-1){27}}

\put (-12,6){$\scriptstyle a_3$}
\put (10,-15){$\scriptstyle a_3$}
\put (-12,80){$\scriptstyle a_5$}
\put (10,100){$\scriptstyle a_5$}
\put (145,80){$\scriptstyle a_4$}
\put (145,6){$\scriptstyle a_2$}
\put (125,-15){$\scriptstyle a_2$}
\put (125,100){$\scriptstyle a_4$}
\end{picture}
}

\begin{equation}\label{monpicture}
\begin{array}{l}
{\scriptstyle a_1 }
\left \{
\ \ \ \
\begin{picture}(150,55)
\put(0,0){\usebox{\build}}
\end{picture}
\right. \\ \\
\hskip .54in {\underbrace{\hskip 2in}_{a_6}}
\end{array}
\end{equation}
Note that the $a_i$ in the diagram measure the length,
and not the number of vertices.  Note also that the minimality, and in
particular the inequalities (\ref{cond for min}), guarantee that these
``cuts'' do not overlap.  In particular, any monomial that is removed
falls unambiguously into exactly one of the removed regions.

Now let us define $I(a_1,a_2,a_3,a_4,a_5,a_6)$ to be the ideal generated by
the remaining monomials, after removing the corners as described above.  We
first make the following claim:

\medskip

\noindent {\bf Claim:} {\em $I(a_1,a_2,a_3,a_4,a_5,a_6) =
(a_1,a_2,a_3,a_4,a_5,a_6)$, that is that $I(a_1,a_2,a_3,a_4,a_5,a_6)$ is
precisely the ideal of the corresponding  tetrahedral curve.}

\medskip

It is clear that the monomial generators that we have removed from
$(a_1,0,0,0,0,a_6)$ are precisely those generators that are not in
$(a_1,a_2,a_3,a_4,a_5,a_6)$.  This proves $\subseteq$.

We now prove the reverse inclusion.  Since $(a_1,a_2,a_3,a_4,a_5,a_6)$ is a
monomial ideal and $(a_1,a_2,a_3,a_4,a_5,a_6) \subset (a_1,0,0,0,0,a_6)$, we
see that the two ideals in the statement of the claim agree in degrees
$\leq a_1 + a_6$.  The danger is that $(a_1,a_2,a_3,a_4,a_5,a_6)$ could have a
monomial minimal generator of larger degree.   We now show that this does not
occur.

Suppose that $M$ were such a generator, of degree $>a_1+a_6$.  $M$
cannot be a multiple of a minimal generator of
$I(a_1,a_2,a_3,a_4,a_5,a_6)$, by the first inclusion.   But on the other hand,
$M$ is contained in $(a_1,0,0,0,0,a_6)$.    Therefore $M$ is a multiple of
one of the monomials that we have removed, say $N$.

We will consider the case where $N$ is in the removed ``Northwest corner;"
the other cases are identical.  We have removed $a_5$ diagonals from this
corner, where $a_5$ is the exponent of the component $(b,d)$.  Clearly if we
multiply $N$ by either $a$ or $c$, we do not produce a monomial that is in
$(a_1,a_2,a_3,a_4,a_5,a_6)$.  (The problem is in the component $(b,d)$.)

More generally, suppose that $N$ lies on the $k$-th diagonal from the
``Northwest corner" ($k \leq a_5$).  Then $N$ has the form
$a^{\ell_1}b^{\ell_2}c^{\ell_3}d^{\ell_4}$ where
\begin{equation}\label{ineqs1}
\begin{array}{rcl}
\ell_1 + \ell_2 & = & a_1, \\
\ell_3 + \ell_4 & = & a_6, \\
\ell_2 + \ell_4 & = & k-1 \ (< a_5)
\end{array}
\end{equation}
Suppose that we multiply $N$ by a monomial $a^{k_1}b^{k_2}c^{k_3}d^{k_4}$.
    The result is, of course,
$a^{k_1+\ell_1}b^{k_2+\ell_2}c^{k_3+\ell_3}d^{k_4+\ell_4}$.  This will be in
$(a_1,a_2,a_3,a_4,a_5,a_6)$ if and only if
\begin{equation}\label{ineqs2}
\begin{array}{rcl}
k_1+\ell_1 + k_2+\ell_2 & \geq & a_1 \\
k_1+\ell_1 + k_3+\ell_3 & \geq & a_2 \\
k_1+\ell_1 + k_4+\ell_4 & \geq & a_3 \\
k_2+\ell_2 + k_3+\ell_3 & \geq & a_4 \\
k_2+\ell_2 + k_4+\ell_4 & \geq & a_5 \\
k_3+\ell_3 + k_4+\ell_4 & \geq & a_6.
\end{array}
\end{equation}

We make the subclaim that (\ref{ineqs2}) holds if and only if
$k_2+k_4 \geq a_5-(k-1)$.
If (\ref{ineqs2}) holds then from (\ref{ineqs2}) and
(\ref{ineqs1}) respectively we have
\[
k_2 + k_4 \geq a_5 - \ell_2 - \ell_4 = a_5 - (k-1)
\]
as desired.  Conversely, assume that $k_2+k_4 \geq a_5-(k-1)$. The first and
last inequalities of (\ref{ineqs2}) are immediate from (\ref{ineqs1}),
without needing our hypothesis. Similarly, none of the other inequalities
apart from the fifth one need our hypothesis.  Consider for instance the
second inequality.  We have
\[
\begin{array}{rcll}
k_1+\ell_1 + k_3+\ell_3 & = & k_1+k_3 + a_1 + a_6 - (k-1) & \hbox{ (from
(\ref{ineqs1}))} \\
& > & k_1+k_3+ a_2+a_4 + a_4+a_5 - (k-1) & \hbox{ (from (\ref{cond for min}))}
\\
& > & k_1+k_3+a_2 + 2 a_4 & \hbox{ (since $k \leq a_5$)} \\
& \geq & a_2.
\end{array}
\]
We have only to show the fifth inequality, and this is immediate using
(\ref{ineqs1}) and our hypothesis.  This completes the proof of our subclaim.

We return to our monomial $M =
a^{k_1+\ell_1}b^{k_2+\ell_2}c^{k_3+\ell_3}d^{k_4+\ell_4}$.  By our subclaim,
$M \in (a_1,a_2,a_3, a_4,a_5,a_6)$ if and only if $k_2+k_4 \geq a_5-(k-1)$.
In order for $M$ to have a chance to be a minimal generator of $(a_1,a_2,a_3,
a_4,a_5,a_6)$, we have to choose $k_1,k_2,k_3,k_4$ as small as possible.
This means that we may assume that $M$ has the form
\[
a^{\ell_1} b^{\ell_2 + k_2} c^{\ell_3} d^{\ell_4+k_4}
\]
where $k_2+k_4 = a_5-(k-1)$.   We have to show that this is in fact a
multiple of an element, $N$, of $I(a_1,a_2,a_3,a_4,a_5,a_6)$.  We assert that
the following is the desired element:
\[
N' = a^{\ell_1-k_2} b^{\ell_2 + k_2} c^{\ell_3-k_4} d^{\ell_4+k_4}.
\]
To see this, consider the matrix (\ref{bigmatrix}). Recall that the monomial
$N = a^{\ell_1}b^{\ell_2}c^{\ell_3}d^{\ell_4}$ lies in the removed Northwest
corner (see diagram  (\ref{monpicture})), on the $k$-th diagonal.  Each step
South represents a decrease of $\ell_1$ by 1 and an increase of $\ell_2$ by
1, while each step East represents a decrease of $\ell_3$ by 1 and an
increase of $\ell_4$ by 1.  The monomial $N'$ (accepting temporarily that the
exponents are non-negative) represents a move from $N$ of $k_2$ steps South
and $k_4$ steps East.  Each such step moves one to the next diagonal.  But
$k_2+k_4 = a_5 -(k-1)$, so the result of the $k_2+k_4$ steps shows that $N'$
lies on the border of (\ref{monpicture}), on the diagonal in the Northwest
corner.  Hence $N' \in I(a_1,a_2,a_3,a_4,a_5,a_6)$ and so also $M \in
I(a_1,a_2,a_3,a_4,a_5,a_6)$.  This completes the proof of the claim.  In
particular, we have shown that $(a_1,a_2,a_3,a_4,a_5,a_6)$ is generated in
degree $a_1+a_6$, and the value of $\beta_1$ is clear from the construction
(see the diagram (\ref{monpicture})).
\smallskip

Having found the minimal generators of our S-minimal curve we can use the
theory of cellular resolutions as developed in \cite{BSt}. Hence, we will first
define an appropriate regular cell complex $X$. Together with a choice of an
incidence function, the cell complex determines a complex $\FF_X$ of free
$R$-modules. Finally, we will show that it is acyclic and resolves the ideal of
our curve.

As above we begin by describing the  cell complex $X_0$   for the curve
$(a_1,0,0,0,0,a_6)$. Its geometric realization will be a rectangle in the
$(i,j)$-plane. To this end we will identify the integer point $(i,j)$ with the
monomial $a^j b^{a_1-j} c^{a_6-i} d^i $.  The vertices of $X_0$ are all the
integer points $(i, j)$ satisfying $0 \leq i \leq a_6$ and $0 \leq j \leq a_1$.
   The edges of $X_0$ are the relative interiors of all ``horizontal'' line
segments with endpoints $(i-1, j)$ and $(i, j)$ and all ``vertical'' line
segments with endpoints $(i, j-1)$ and $(i, j)$ where all the endpoints are
vertices of $X_0$. The $2$-dimensional faces of
$X_0$ are the relative interiors of all squares with vertices $(i-1, j-1), (i,
j - 1), (i, j), (i-1, j)$, all being vertices of $X_0$. Together with the empty
set, $X_0$ is clearly a finite regular cell complex (cf., e.g., \cite{BH},
Section 6.2 for the definition). The empty set has dimension $-1$.

The cell complex $X := X(a_1,a_2,a_3,a_4,a_5,a_6)$ of the curve
$(a_1,a_2,a_3,a_4,a_5,a_6)$ is a sub-complex of $X_0$ obtained by deleting
corners of $X_0$. More precisely, the vertices of $X$ are the points $(i, j)$
satisfying
\begin{equation} \label{eq-gen}
\begin{array}{ccccc}
a_3 & \leq & i + j & \leq & a_1 + a_6 - a_4\\
a_5 - a_1 & \leq & i - j & \leq & a_1 - a_2 \\
\end{array}
\end{equation}
The edges and $2$-dimensional faces of $X$ are by definition the
corresponding  faces of $X_0$ whose edges all belong to $X$.

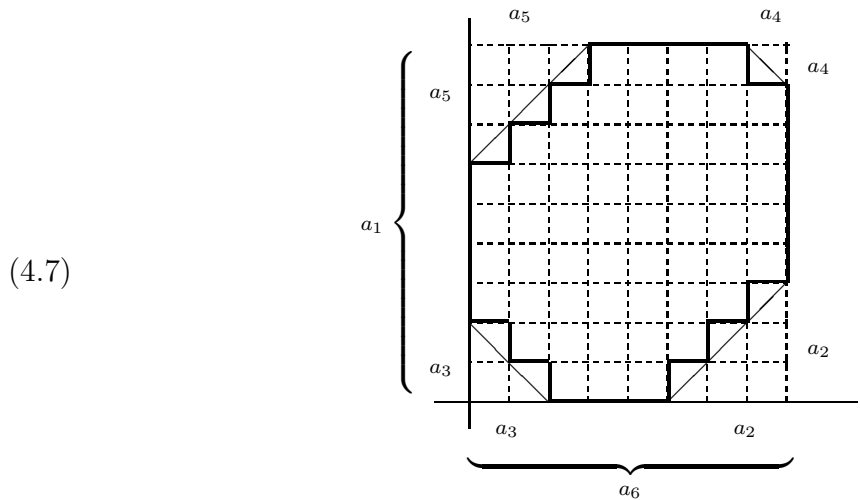
\newsavebox{\cell}
\savebox{\cell}(200,100)[tl]
{
\begin{picture}(200,200)
\put (0,-10){\line (0,1){155}}

\put (-13,0){\line (1,0){160}}

\put (15,0){\dashbox{2}(0,135)[t]}
\put (30,0){\dashbox{2}(0,135)[t]}
\put (45,0){\dashbox{2}(0,135)[t]}
\put (60,0){\dashbox{2}(0,135)[t]}
\put (75,0){\dashbox{2}(0,135)[t]}
\put (90,0){\dashbox{2}(0,135)[t]}
\put (105,0){\dashbox{2}(0,135)[t]}
\put (120,0){\dashbox{2}(0,135)[t]}

\put (0,135){\dashbox{2}(120,0)[t]}
\put (0,120){\dashbox{2}(120,0)[t]}
\put (0,105){\dashbox{2}(120,0)[t]}
\put (0,90){\dashbox{2}(120,0)[t]}
\put (0,75){\dashbox{2}(120,0)[t]}
\put (0,60){\dashbox{2}(120,0)[t]}
\put (0,45){\dashbox{2}(120,0)[t]}
\put (0,30){\dashbox{2}(120,0)[t]}
\put (0,15){\dashbox{2}(120,0)[t]}

\put (0,90){\line (1,1){45}}
\put (0,30){\line (1,-1){30}}
\put (105,135){\line (1,-1){15}}
\put (75,0){\line (1,1){45}}

\put (0,30){\rule{1pt}{60pt}}
\put (0,90){\rule{15pt}{1pt}}
\put (15,90){\rule{1pt}{15pt}}
\put (15,105){\rule{15pt}{1pt}}
\put (30,105){\rule{1pt}{15pt}}
\put (30,120){\rule{15pt}{1pt}}
\put (45,120){\rule{1pt}{15pt}}
\put (45,135){\rule{60pt}{1pt}}
\put (105,120){\rule{1pt}{15pt}}
\put (105,120){\rule{15pt}{1pt}}
\put (120,45){\rule{1pt}{75pt}}
\put (105,45){\rule{15pt}{1pt}}
\put (105,30){\rule{1pt}{15pt}}
\put (90,30){\rule{15pt}{1pt}}
\put (90,15){\rule{1pt}{15pt}}
\put (75,15){\rule{15pt}{1pt}}
\put (75,0){\rule{1pt}{15pt}}
\put (30,0){\rule{45pt}{1pt}}
\put (30,0){\rule{1pt}{15pt}}
\put (15,15){\rule{15pt}{1pt}}
\put (15,15){\rule{1pt}{15pt}}
\put (0,30){\rule{15pt}{1pt}}

\put (-15,11){$\scriptstyle a_3$}
\put (10,-12){$\scriptstyle a_3$}
\put (-15,115){$\scriptstyle a_5$}
\put (15,145){$\scriptstyle a_5$}
\put (110,145){$\scriptstyle a_4$}
\put (128,18){$\scriptstyle a_2$}
\put (100,-12){$\scriptstyle a_2$}
\put (128,125){$\scriptstyle a_4$}
\end{picture}
}

\begin{equation}\label{monpicture2}
\begin{array}{l}
{\scriptstyle a_1 }
\left \{
\ \ \ \
\begin{picture}(150,70)
\put(0,35){\usebox{\cell}}
\end{picture}
\right. \\ \\
\hskip .56in {\underbrace{\hskip 1.7in}_{a_6}}
\end{array}
\end{equation}

In the above diagram, the full rectangle indicates the cell complex $X_0$,
while the bold face segments are the boundary of the cell complex $X$.
Again it is clear that
$X$ is a finite regular cell complex. For a face
$F$ of $X$ we set $m_F$ the least common multiple of the vertices of $F$. The
exponent vector of the monomial $m_F$ in $\ZZ^4$ is called the degree of the
face $F$.

Now we fix an incidence function $\eps(F,F')$ on pairs of faces of $X$ (see for
instance \cite{BSt}). Then the  cellular complex $\FF_X$ is the $\ZZ$-graded
complex
\begin{equation*}
\begin{CD}
\FF_X: \quad \quad \quad 0 @>>> F_{2}  @>{\ffi_{2}}>> F_1  @>{\ffi_1}>> F_0
@>{\ffi_{0}}>>0
\end{CD}
\end{equation*}
where
$$
F_i = \bigoplus_{F \in X, \dim F = i} R(- \deg m_F) = \bigoplus_{F \in X, \dim
F = i} R e_F
$$
and the differential is given by
\begin{equation} \label{eq-diff}
\ffi_i(F) = \sum_{F' \in X, \dim F' = i-1} \eps(F, F') \frac{m_F}{m_{F'}} e_F.
\end{equation}
By the construction of this complex, the image of $\ffi_0$ is the
ideal
$(a_1,a_2,a_3,a_4,a_5,a_6)$, thanks to the computation above of the minimal
generators. We want to show that
$\FF_X$ provides a minimal free resolution of this ideal. To this end denote
for $b = (b_0,\ldots,b_3) \in
\ZZ^4$ by $X_{\leq b}$ the sub-complex of $X$ on the vertices of degree $\leq
b$. Here, we use the partial ordering on $\ZZ^3$ induced by comparing
componentwise.  It is easy to see that for all degrees $b \in \ZZ^4$ the 
geometric realization of the  cell complex $X_{\leq b}$ is contractible,
thus $X_{\leq b}$ is acyclic. Hence,
\cite{BSt}, Proposition 1.2, shows that $\FF_X$ is a free resolution of the
ideal $(a_1,a_2,a_3,a_4,a_5,a_6)$.

Next, we observe  for every $i$-dimensional face $F$ of $X$ that $m_F$ has
degree $a_1 + a_6 + i$. It follows immediately that $\FF_X$ is a minimal free
resolution. Moreover, it is easy to see that
$X$ contains $(a_1 + 1) (a_6 + 1) - \sum_{i=2}^5 \frac{a_i (a_i + 1)}{2}$
vertices,
$2 a_1 a_6 + a_1 + a_6  - \sum_{i=2}^5 a_i (a_i + 1)$ edges, and $a_1 a_6  -
\sum_{i=2}^5 \frac{a_i (a_i + 1)}{2}$ $2$-dimensional faces. Hence the complex
$\FF_X$ giving the free resolution of $C$ is of the form as claimed.
\end{proof}

Combining this result with the Auslander-Buchsbaum formula we get:

\begin{corollary} \label{cor-acm}
An $S$-minimal tetrahedral curve is arithmetically Cohen-Macaulay if and only
if it is trivial.
\end{corollary}

\begin{remark} \label{comparison}
(i)
In \cite{BSt} Bayer and Sturmfels introduce the hull resolution. They show that
the hull resolution of a monomial ideal is a free resolution which is not
necessarily minimal. The hull resolution contains as sub-complex the so-called
Scarf complex which is easy to compute. For generic monomial ideals the Scarf
complex and the hull resolution agree and provide a minimal free resolution.
However, the defining ideal of a tetrahedral curve is typically not generic (in
the sense of Bayer and Sturmfels). For example, the ideal $(1,0,0,0,0,1)$
corresponding to a pair of skew lines has a Scarf complex of length $2$ while
its minimal free resolution has length $3$ (cf.\ \cite{Ta}, Example 3.4.20).
Tappe has also computed the hull resolution of the ideals $(a_1,0,0,0,0,a_6)$
for $a_1, a_6 \leq 3$. In these cases the hull resolution does give the minimal
free resolution. However, it is not clear (to the authors) if this is true for
all $S$-minimal curves.

Note that the computation of the hull resolution requires one to compute the
convex  hull of as many points as the ideal has minimal generators. This is a
rather non-trivial task if the number of points is large. On the other hand our
description of the minimal free resolution of $S$-minimal curves uses a very
simple cell complex which gives a more direct approach.

(ii) After we had found the cell complex that determines the minimal free
resolution of an $S$-minimal curve we realized that Schwartau \cite{Sw} uses a
similar description in order to compute the minimal free resolution of any
tetrahedral curve with $a_2 = a_5 = 0$. However, much of the theory that we
have used here was not available during the writing of \cite{Sw}, and we found
that proof to be less ``clean'' even in the restricted setting.
\end{remark}

As another consequence of Theorem \ref{thm-res} we get some information about
the Hartshorne-Rao module of any tetrahedral curve. Recall that the
Hartshorne-Rao module of a curve $C$ is the graded module
$$
M(C) := \oplus_{j \in \ZZ} H^1(\cI_C(j)).
$$

\begin{corollary} \label{cor-HR}
The $K$-dual of the Hartshorne-Rao module of a tetrahedral curve is generated
in one degree.
\end{corollary}

\begin{proof} Since the Hartshorne-Rao module is up to degree shift invariant
in an even liaison class, it suffices to show the claim for a non-trivial
$S$-minimal curve $C$.  By local duality we have the following graded
isomorphism
$$
M(C)^{\vee} \cong \Ext^3(R/I_C, R) (-4).
$$
The latter module can be computed from the minimal free resolution of $I_C$.
Since it is linear the claim follows.
\end{proof}

\begin{remark} \label{rem-HR}
(i) The conclusion of the last result is in general not true for the \linebreak
Hartshorne-Rao module itself. Indeed, the Hartshorne-Rao module of the curve
\linebreak
$(5,2,2,1,1,5)$ has minimal generators in two different degrees, namely 4 and
5.

(ii) The curve in (i) also shows that in general we cannot link a tetrahedral
curve in an odd number of steps to a tetrahedral curve.
\end{remark}

  From a cohomological point of view the simplest curves that are not
arithmetically Cohen-Macaulay are the arithmetically Buchsbaum curves. The
tetrahedral curves among them have special properties.
\begin{corollary} \label{cor-HR-abm}
A tetrahedral curve is arithmetically Buchsbaum if and only if its \linebreak
Hartshorne-Rao module satisfies
$$
M(C) \cong K^m(t) \quad \mbox{for some integers $m,t$ where $m \geq 0$}.
$$
\end{corollary}

\begin{proof}
This follows immediately by Corollary \ref{cor-HR}.
\end{proof}

\begin{remark} \label{schwartau knew}
The fact that such arithmetically Buchsbaum curves exist was known already to
Schwartau, who constructed them with liaison addition.  The fact that they are
the {\em only} arithmetically Buchsbaum curves among the tetrahedral curves is
new.  
\end{remark}


\section{Minimality in the Even Liaison Class, and Applications}

Now we want to show that for tetrahedral curves the concept of $S$-minimality
and minimality in its even liaison class agree.   This is the second main
result of this paper.

\begin{theorem} \label{prop-min}
Let $C$ be a tetrahedral curve which is not arithmetically Cohen-Macaulay. Then
$C$ is $S$-minimal if and only if it is minimal in its even liaison class.
\end{theorem}

\begin{proof}
It suffices to show that every $S$-minimal curve is minimal in its even liaison
class. Denote by $s(C)$ the initial degree of $C$,
$$
s(C) := \min\{j \in \ZZ \s [I_C]_j \neq 0\},
$$
and by $e(C)$ its index of speciality,
$$
e(C) := \max\{j \in \ZZ \s H^2(\cI_C(j)) \neq 0\}.
$$
A result of Lazarsfeld and Rao in \cite{LR} says that $C$ is minimal if
$$
s(C) \geq e(C) + 4.
$$
We will show that this criterion applies to our tetrahedral curves.
We may again assume that $a_6 = \max\{a_i \} > 0$.  The minimal free resolution
of $C$ is
\[
\begin{array}{cccccccccccccccccccccc}
0 & \rightarrow & R^{\beta_3}(-a_1-a_6-2) & \stackrel{\phi_2}{\rightarrow} &
R^{\beta_2}(-a_1-a_6-1) & \rightarrow & R^{\beta_1}(-a_1-a_6) & \rightarrow &
I_C & \rightarrow & 0 \\
&&&&\hfill \searrow && \nearrow \hfill \\
&&&&& K \\
&&&& \hfill \nearrow && \searrow \hfill\\
&&&& \hskip 1.6cm 0 && 0 \hskip 1.cm
\end{array}
\]
where $K$ splits the resolution into two short exact sequences.  We know that
$s(C) = a_1+a_6$, so we have to show that
$h^2({\mathcal I}_C( a_1+a_6-3)) = 0$.  (See Remark \ref{rem-mini} (ii).)
Letting
$\mathcal K$ be the sheafification of
$K$, it is enough to show that $h^3({\mathcal K}(a_1+a_6-3)) = 0$.
The leftmost short exact sequence gives
\[
0 \rightarrow H^2({\mathcal K}(a_1+a_6-3)) \rightarrow H^3({\mathcal
O}^{\beta_3}_{\proj{3}}(-5)) \rightarrow H^3({\mathcal
O}_{\proj{3}}^{\beta_2}(-4)) \rightarrow H^3({\mathcal K}(a_1+a_6-3))
\rightarrow 0.
\]
Since the rightmost short exact sequence gives $h^1({\mathcal I}_C(t)) =
h^2({\mathcal K}(t))$ for all $t$, the above long exact sequence gives
\begin{equation}\label{fromles}
h^3({\mathcal K}(a_1+a_6-3)) = \beta_2 - 4 \beta_3 + h^1({\mathcal I}_C
(a_1+a_6-3)).
\end{equation}
Now, recall from \cite{rao}, Theorem 2.5, that if $C$ is a curve and if $M(C)$
has minimal free resolution
\[
0 \rightarrow L_4 \stackrel{\sigma_4}{\longrightarrow} L_3 \rightarrow L_2
\rightarrow L_1 \rightarrow L_0 \rightarrow M(C) \rightarrow 0
\]
then $I_C$ has minimal free resolution
\[
0 \rightarrow L_4 \stackrel{(\sigma_4,0)}{\longrightarrow} L_3 \oplus F
\rightarrow F_1 \rightarrow I_C \rightarrow 0.
\]
In our case we know the minimal free resolution of $I_C$.  We first claim that
$F = 0$, i.e.\ the last two free modules (and the map between them) in the
minimal free resolution of $I_C$ exactly coincide with the corresponding ones
for $M(C)$.  In the proof of Theorem \ref{thm-res} we have described the maps
in the minimal free resolution of $C$. Denote by $M$ the matrix describing the
map $\phi_2$ after choosing canonical bases for the free modules
$F_i$, where we
think of the columns of $M$ as second syzygies.  We want to show that $M$ does
not have a row of zeros.  But if it did,
this means that there is a first syzygy which does not ``contribute'' to any
second syzygy. This gives a contradiction since the first syzygies correspond
to edges of the cell  complex  $X$, but every edge of $X$ is in the boundary of
a facet of $X$ that corresponds to a second syzygy. (This is a consequence of
the inequalities (\ref{cond for min}).)

Thus we know the end of the minimal free resolution for $M(C)$ has the form
\[
0
\rightarrow R^{\beta_3}(-a_1-a_6-2) \stackrel{\phi_2}{\longrightarrow}
R^{\beta_2}(-a_1-a_6-1)
\rightarrow \dots,
\]
   and as a
result the minimal free resolution for the
dual module (over the field $K$)
$M(C)^\vee$ has the form
\[
\dots \rightarrow {\mathbb F} \rightarrow R^{\beta_2} (a_1+a_6-3)
\stackrel{\phi_2^\vee}{\longrightarrow} R^{\beta_3}(a_1+a_6-2)
\rightarrow M(C)^\vee \rightarrow 0.
\]
It is a basic property of minimal free resolutions that $\mathbb F$ has no
summand of the form
$R(a_1+a_6-3)$. Hence $\dim
M(C)^\vee_{-a_1-a_6+3} = 4 \beta_3 -
\beta_2$, so (\ref{fromles}) gives $h^3({\mathcal K}(a_1+a_6-3)) = 0$ as
desired.
\end{proof}

\begin{corollary} \label{deg-genus}
Assume that $a_6 = \max \{ a_1,\dots,a_6\}$.  Then the degree of an S-minimal
curve of type
$(a_1,\dots,a_6)$ is
$\sum_{i=1}^6 \binom{a_i + 1}{2}$, and the arithmetic genus is
\[
\left [ \sum_{i=1}^6 \binom{a_i +1}{2} \right ](a_1+a_6-1) +1- \binom{a_1 +
a_6 +2}{3}.
\]
\end{corollary}

\begin{proof}
The degree statement is obvious.  As for the arithmetic genus, let us denote
it by $g$.  It follows easily from the calculations in the proof of Theorem
\ref{prop-min} that
\[
h^0(\mathcal I_C (a_1+a_6-1)) = h^1(\mathcal
I_C(a_1+a_6-1)) = h^2(\mathcal I_C(a_1+a_2-1)) = 0.
\]
  From the cohomology of the exact sequence
\[
0 \rightarrow \mathcal I_C \rightarrow \mathcal O_{\mathbb P^3} \rightarrow
\mathcal O_C \rightarrow 0
\]
(twisted by $a_1+a_6-1$) it follows that
\[
h^0(\mathcal O_C (a_1+a_6-1)) = \binom{a_1+a_6+2}{3}.
\]
We also  have $h^1(\mathcal O_C(t)) = h^2(\mathcal I_C(t))$ for all
$t$.  But the Riemann-Roch theorem gives that
\[
h^0(\mathcal O_C(a_1+a_6-1)) =
(\deg C)(a_1+a_6-1) -g+1 + h^1(\mathcal O_C (a_1+a_6-1)).
\]
Combining the above calculations gives the result.
\end{proof}

The following  theorem can be viewed as a clarification of Schwartau's theorem
(see Theorem \ref{Sw} and Remark \ref{what phil meant}) and its proof
is new and
uses the methods of our paper.  The main tool is Lemma \ref{lem-max-fac}, but
we have put it after Theorem \ref{prop-min} because we use the fact that an
S-minimal curve is not ACM, and this follows from the fact that a curve that is
minimal in its even liaison class is not ACM.

\begin{theorem}[invariant Schwartau] \label{thm-corr-sw}
The ideal
\[
(a,b)^{a_1} \cap (a,c)^{0} \cap (a,d)^{a_3} \cap (b,c)^{a_4} \cap
(b,d)^{0} \cap (c,d)^{a_6}
\]
defines an \acm curve in $\proj{3}$ if and only if
\begin{itemize}
\item[\underline{Case 1}.]  $a_1,a_3,a_4,a_6 >0$: \quad $a_1+a_6 = a_3+a_4+
\epsilon, \hbox{ where } \epsilon \in \{ -1,0,1 \}$. \\
\item[\underline{Case 2}.]  Exactly one of $a_1,a_3,a_4,a_6$ is zero, say
$a_i$: \quad $ a_{7-i}+1 \geq $ the sum of the weights of the lines meeting the
line $l_i$.
\\
\item[\underline{Case 3}.]  At least two of $a_1,a_3,a_4,a_6$ are zero: \quad
the curve is connected.
\end{itemize}
\end{theorem}

\medskip

\begin{proof}
    For convenience we will call
a curve of the form given in Theorem \ref{thm-corr-sw} a {\em Schwartau
curve.}  If
$C$ is a Schwartau curve then its components can be represented by a square:
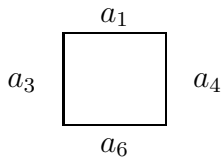
\[
\begin{array}{rcl}
&a_1 \\
a_3 &\fbox{\rule[-4mm]{0cm}{1cm}\phantom{sssssss}} & a_4 \\
& a_6
\end{array}
\]
If three of the integers $a_i$ are zero then the ideal is a power
of a complete intersection, and is automatically ACM and connected.  Suppose
that two of the integers $a_i$ are zero.  If $C$ is ACM then it is a general
fact that it must be connected.  Conversely, suppose that $C$ is connected.
This means that the two non-zero values of $a_i$ represent adjacent
sides of the
square.  Hence the ideal involves only three variables.  The fourth variable
is thus a non zero-divisor, and reducing modulo this variable gives the
saturated ideal of a zeroscheme in $\proj{2}$.  Hence in this case the
Schwartau curve is ACM.  This proves Case 3.

Now suppose that one of the integers is zero, say $a_i$.  (Note that Schwartau
took $a_3 = 0$.)  As long as
$C$ is not S-minimal, we can reduce any maximal facet, by Lemma
\ref{lem-max-fac}.  But with $a_i = 0$, there remain only two facets, whose
weights are the sums of two consecutive sides of the square.  (For example, if
$a_3 = 0$ then the weights of the facets are
$a_4+a_6$ and $a_1+a_4$.)

If $C$ is ACM then no reduction will ever be S-minimal.  We repeatedly reduce
the facet of maximal weight until we obtain another zero for the weight of an
edge, and the fact that $C$ is ACM means that after obtaining this zero, the
remaining two edges must be connected (by Case 3 above).  This means that the
``middle'' of the three non-zero edges, $a_{7-i}$, must not be the first to
reach 0 unless there is a tie.  But this middle edge is involved in every
reduction!  Note that since we are always reducing the facet of maximal weight,
we eventually reach the point where the weights of the non-middle edges differ
by at most 1.  It follows that $a_{7-i} +1 \geq $ the sum of the weights of the
two ``non-middle'' edges, as claimed.

Conversely, suppose that $a_{7-i} +1 \geq $ the sum of the weights of the
two ``non-middle'' edges.  Then clearly $a_{7-i}$ is the edge of maximal
weight.  The condition (b) of Lemma \ref{lem-max-fac} is trivial to check, so
$C$ is not minimal and we can reduce the facet of maximal weight.  In doing so,
both sides of the inequality in Case 2 are reduced by 1, so the inequality
still holds for the new curve.  Again, reducing the facet of maximal weight
each time eventually leads to a balance in the two non-middle edges, and the
given inequality guarantees that the curve resulting when a second 0 is
obtained, will be connected, hence ACM.  So $C$ is ACM.

We now consider the case where none of the $a_i$ are 0.  Again without loss of
generality assume that $a_6 = \max \{ a_1, a_3, a_4,a_6 \}$.  One can
check that
then $w = \max \{ a_3+a_6, \linebreak a_4+a_6 \}$.   It follows that when we
reduce a facet of maximal weight (if it is possible), we simultaneously reduce
$\max\{a_1,a_6 \}$ and $\max \{ a_3,a_4\}$ by 1.  In particular, we
simultaneously reduce $(a_1+a_6)$ by 1 and $(a_3 +a_4)$ by 1.

   If $C$ is ACM then it is not S-minimal, so by Lemma \ref{lem-max-fac} we have
\[
a_1+a_6 \leq \max \{ a_3+a_6, a_4+a_6 \}, \hbox{ i.e.\ $a_1 \leq \max \{
a_3,a_4\}$}.
\]
    We conclude that if $C$ is ACM then the largest two weights are on
consecutive edges.  That is, the facet of maximal weight is one of the two
facets involving $a_6$, and we can reduce that facet.  Since $C$ is ACM, we
will never reach an S-minimal curve in this process, and the result of
performing a reduction will again be an ACM curve.   Hence after each step in
the reduction, the two edges of largest weight are consecutive and each can be
reduced by 1 by a new reduction.

Suppose that we have reached the point where the application of Proposition
\ref{reduce} reduces an $a_i$ to zero.  If two edges are reduced to 0
simultaneously, the remaining two edges are consecutive and have weight 1
each (at this point).  So because of our procedure, it follows that $a_1+a_6 =
a_3+a_4$.

If one edge is reduced to 0, suppose that it is supported on $\ell_i$.  Then at
this point $\ell_{7-i}$ has weight 1, and the resulting curve is ACM.  But
$\ell_{7-i}$ is then the middle edge, so by Case 2 we have the other two edges
having weight 1 each.  It follows that $a_1+a_6 = a_3+a_4 + \epsilon$ as
claimed.

Conversely, assume that $a_1+a_6 = a_3+a_4+ \epsilon$ for $ \epsilon =
-1,0,1$.  We want to show that a reduction can be performed, i.e.\ $C$
is not minimal.  Without loss of generality we still assume that $a_6$ is the
largest weight.  We claim that $a_1 \leq \max \{ a_3, a_4 \}$.  Indeed, if
$a_6 \geq a_1 > \max \{ a_3,a_4\}$ then $a_1+a_6 \geq a_3+a_4+2$, contradicting
the hypothesis.  But then if $w$ is the maximal weight of a facet, we have
$w = a_6 + \max \{ a_3,a_4\}$ and $a_1+a_6 \leq w$.  Therefore, by Lemma
\ref{lem-max-fac}, a reduction can be performed.

Hence our numerical assumption guarantees that at least as long as the four
integers are positive, we can always perform a reduction, reducing the larger
of $a_1$ and $a_6$ by 1 and the larger of $a_3$ and $a_4$ by 1.  Suppose that
$a_3$ is the first to be reduced to zero.  Then at that point $a_4=1$ and
$a_1+a_6$ is either 1 or 2.  Either way, we are in the previously studied
cases and $C$ is ACM.
\end{proof}

We can also characterize the $S$-minimal arithmetically Buchsbaum curves.

\begin{corollary} \label{cor-abm}
Let  $C$ be an $S$-minimal tetrahedral curve  defined by
$(a_1,a_2,a_3,a_4,a_5,a_6)$. Then $C$ is arithmetically Buchsbaum and not
arithmetically Cohen-Macaulay if and only if there are integers $i \neq j$ in
$\{1, 2, 3\}$ such that
$$
a_i + 1 = a_j = a_{7-j} = a_{7-i} + 1
$$
and the two remaining weights are zero.
\end{corollary}

\begin{proof}
According to Corollary \ref{cor-acm} we may assume that $a_6 = \max \{a_i\} >
0$. We will use the notation of Theorem \ref{thm-res} and its proof. Since $C$
is minimal, a result of Martin-Deschamps  and Perrin says that  the transpose
of the matrix describing the map $\ffi_2$ is a minimal presentation matrix of
$\Ext^3 (R/I, R)$.  (See \cite{rao} Theorem 2.5 and \cite{MDP} Proposition
IV.4.4.)  We also saw this directly, in the case of tetrahedral curves, in the
proof above. By Corollary
\ref{cor-HR-abm},
$C$ is arithmetically Buchsbaum if and only if the latter module is a (shifted)
direct sum of copies of $K$. Using the cell complex $X$ that governs the
minimal free resolution of
$C$, this is equivalent to the fact that every edge of $X$ belongs to the
boundary of exactly one facet. Now the claim follows from the description of
$X$.  Indeed, after removing diagonals from the corners of (\ref{bigmatrix})
and translating to the notation of the cell complex $X$, we must be left with a
diagonal of facets and no additional edges, and this can only happen if we only
remove diagonals from two opposite corners of (\ref{bigmatrix}) in the
prescribed way.
\end{proof}

\begin{remark} \label{rem-mini}
(i) The arguments in the proof of Theorem \ref{prop-min}, coupled with a recent
result of Strano  \cite{Strano}, Theorem 1, show the following fact:

If a non-arithmetically Cohen-Macaulay curve $C \subset \PP^3$ has a linear
resolution then it is minimal in its even liaison class if and only of $s(C)
\geq e(C) + 4$.

(ii) There are tetrahedral curves that have a linear resolution, but are not
minimal in their even liaison class. An example is given by the curve with
weights $(5, 1, 3, 2, 2, 5)$.  This shows that our task in the proof of Theorem
\ref{prop-min}, of showing that $h^2({\mathcal I}_C (a_1+a_6-3)) = 0$, cannot
be shown merely by invoking the linearity of the resolution.  Indeed, this only
guarantees
$h^2({\mathcal I}_C (a_1+a_6-2)) = 0$.

(iii) Sometimes a non-trivial $S$-minimal curve is the unique minimal
curve in its even liaison class. This is true if $s(C) \geq e(C) + 5$ 
(thanks to
\cite{LR}). Examples of such curves have weights $(m, 0, 0, 0, 0, k)$ 
where $m, k
\geq 2$.
\end{remark}


\section{Unobstructedness of some curves} \label{sec-unobstr}

We will show that some of the minimal tetrahedral curves correspond to smooth
points of a component of the corresponding Hilbert scheme that has the
``expected'' dimension.

One starting point is the following result whose proof is an adaptation of
Dolcetti's method in \cite{Dolcetti}. Note that if $d_1$ is the degree of the
first non-zero component of a graded module, $M$, of finite length and $d_2$ is
the degree of the last non-zero component of $M$, then the {\em
diameter} of $M$
is $d_2-d_1+1$.

\begin{proposition} \label{lem-dolcet}
Let $C \subset \PP^3$ be a curve with a linear resolution. If the
diameter of its
Hartshorne-Rao module is at most two then its normal sheaf $\cN_C$ satisfies
$$
H^1(\cN_C) = 0.
$$
Therefore, the corresponding component of the Hilbert scheme is
generically smooth
of dimension $4 \cdot \deg C$.
\end{proposition}

\begin{proof}
By assumption the minimal free resolution of $C$ is of the form
\begin{equation*}
   0 \to R^{\beta_3}(-s-2)  \to R^{\beta_2}(-s-1)  \to
R^{\beta_1}(-s) \to I_C \to 0.
\end{equation*}
It implies
$$
H^1(\cI_C(s-2)) \neq 0 \quad \mbox{and} \quad H^1(\cI_C(j)) = 0 \;\;
\mif j  \geq
s-1
$$
as well as
$$
H^2(\cI_C(j)) = 0 \quad \mbox{for all} \; j \geq s-2.
$$
Thus, the assumption on the diameter provides
$$
H^1(\cI_C(j)) = 0 \quad \mbox{for all} \; j \leq s-4.
$$
Sheafifying the resolution above we define the vector bundle $\cE$ as follows
$$
\begin{array}{ccccc} \label{eq-E}
   0 \to \cOP^{\beta_3}(-s-2)  \to & \cOP^{\beta_2}(-s-1) & \longrightarrow &
\cOP^{\beta_1}(-s) & \to \cI_C \to 0 \\
&\hfill \searrow && \nearrow \hfill \\
&& \cE & \\
& \hfill \nearrow && \searrow \hfill\\
& 0 && 0
\end{array}
$$
Thus, we have
$$
H^1_*(\cE) = 0 \quad \mbox{and} \quad H^2_*(\cE) = H^1_*(\cI_C).
$$
Furthermore, dualizing the exact sequence on the left-hand side and
then tensoring
by $\cE$ we get the exact sequence
$$
0 \to \cE^* \otimes \cE \to \cE^{\beta_2}(s+1) \to \cE^{\beta_3}(s+2) \to 0.
$$
Taking cohomology we obtain
$$
H^2(\cE^* \otimes \cE) = 0.
$$
Now, the cohomology of the exact sequence
$$
0 \to \cE^* \otimes \cE \to (\cE^*)^{\beta_1}(-s) \to \cE^* \otimes \cI_C \to 0
$$
provides
$$
(H^1(\cE^*(-s)))^{\beta_1} \to H^1(\cE^* \otimes \cI_C) \to H^2(\cE^* \otimes
\cE) = 0.
$$
Since we have by duality
$$
H^1(\cE^*(-s)) \cong [\Ext(R/I_C, R)]_{-s} \cong H^1(\cI_C(s-4)) = 0
$$
we conclude
$$
H^1(\cE^* \otimes \cI_C) = 0.
$$
The sequence defining $\cE$ also gives the exact sequence
$$
\Ext^1(\cE, \cI_C) \to \Ext^2(\cI_C, \cI_C) \to \Ext^2(\cOP^{\beta_1}(-s),
\cI_C).
$$
Using
$$
\Ext^2(\cOP^{\beta_1}(-s), \cI_C) \cong (H^2(\cI_C(s)))^{\beta_1} = 0
$$
and
$$
\Ext^1(\cE, \cI_C) \cong H^1(\cE^* \otimes \cI_C) = 0
$$
we get
$$
0 = \Ext^2(\cI_C, \cI_C) \cong H^1(\cN_C),
$$
as claimed.
\end{proof}

Since minimal tetrahedral curves have a linear resolution we would like to know
the ones whose Hartshorne-Rao module has diameter at most two. The curves with
diameter one have already been characterized in Corollary \ref{cor-abm}. The
diameter two case is described below:

\begin{lemma} \label{lem-diam-two}
Let $C$ be a minimal tetrahedral curve. Then its Hartshorne-Rao module has
diameter two if and only if $C$ is isomorphic to one of the curves
\begin{align*}
&(k, k-1, 0, 0, k-1, k+1) \quad \mbox{where} \quad k \geq 1, \quad \mbox{or} \\
&(k, k-2, 0, 0, k-1, k) \quad \mbox{where} \quad k \geq 2.
\end{align*}
\end{lemma}

\begin{proof}
Let $C$ be the minimal curve
defined by the tuple $(a_1,a_2,a_3,a_4,a_5,a_6)$.  Without loss of
generality assume that $a_6$ is the largest entry. We will use the
description of
the minimal free resolution of $C$ given in Theorem \ref{thm-res}. Denote by
$\psi$ the dual of the last map in this resolution, i.e.\ $\psi =
\ffi_2^*$. Thus,
we have the minimal presentation
\begin{equation*}
\begin{CD}
   @>>> R^{\beta_3}(a_1+a_6+1)  @>{\psi}>> R^{\beta_2}(a_1+a_6+2)  @>>>
\Ext^3(R/I_C, R) @>>> 0
\end{CD}
\end{equation*}
Since $\Ext^3(R/I_C, R)$ is the $K$-dual of the Hartshorne-Rao module
$M_C$ of $C$,
the diameter of the latter is at most two if and only if
\begin{equation} \label{eq-diam}
[\Ext^3(R/I_C, R)]_{-a_1-a_6}  =  0.
\end{equation}
Denote by $X := X(a_1,a_2,a_3,a_4,a_5,a_6)$ the cell complex of the curve $C$.
Recall from the proof of Theorem \ref{thm-res} that the facets of $X$
correspond
to minimal generators of $\Ext^3(R/I_C, R)$ and that $\im \psi$ has a system of
minimal generators consisting of binomials and monomials (cf.\ formula
(\ref{eq-diff})). These generators
correspond to the edges of $X$.  Below, we will always refer to this system of
generators. Moreover, we will use Diagram (\ref{monpicture2}).

We begin with deriving necessary conditions for Condition (\ref{eq-diam}) being
true. In order to rule out  curves we show two claims.

{\it Claim 1}: If $X$ contains ``3 facets in a row'' then the
diameter of $M_C$ is
greater than two.

Here the assumption means that $X$ contains a subcomplex that is isomorphic to
\linebreak
$X(3,0,0,0,0,1)$ (a ``vertical row'') or $X(1,0,0,0,0,3)$ (a
``horizontal row'').

Now, assume that $X$ contains 3 facets in a horizontal row. Denote by
$e_1, e_2,
e_3$ these facets from left to right. Then, the only minimal generators of $\im
\psi$ involving $c e_1, c e_2, d e_2$, or $d e_3$ are (up to sign)
$$
-c e_1 + d e_2 \quad \mbox{and} \quad -c e_2 + d e_3.
$$
It follows easily that all the monomials $c^2 e_1, cd e_2, d^2 e_3$
do not belong
to $\im \psi$. Thus, condition (\ref{eq-diam}) is not satisfied.

If $X$ contains 3 facets in a vertical row we argue similarly and Claim 1 is
shown.

{\it Claim 2}: If $X$ contains a ``square of 4 facets'' then the diameter of
$M_C$ is greater than two.

More precisely, the assumption means that $X$  contains a subcomplex that is
isomorphic to $X(2,0,0,0,0,2)$. Enumerate its facets counterclockwise by
$e_1,\ldots,e_4$ beginning with $e_1$ in the ``Southwest'' corner.
Then, the only
minimal generators of $\im \psi$ involving $b e_1, c e_1, b e_2, d
e_2, a e_3, d
e_3, a e_4$, or $c e_4$ are (up to sign)
\begin{align*}
& -c e_1 + d e_2, \quad  -b e_2 + a e_3 \\
& -c e_4 + d e_3, \quad  -b e_1 + a e_4.
\end{align*}
Hence, all the monomials $bc e_1, bd e_2, ad e_3, ac e_4$ are not in
$\im \psi$.
Thus, condition (\ref{eq-diam}) is not satisfied and Claim 2 is established.
\\[.5ex]

As the next step we look for the minimal tetrahedral curves whose cell complex
contains neither 3 facets in a row  nor a square of 4 facets.
Let $C$ be
such a curve. Looking at the boundaries of its cell complex $X$ we get in
conjunction with Corollary \ref{s-min}
\begin{align*}
&1 \leq a_6 - a_4 - a_5 \leq 2 \\
&1 \leq a_6 - a_2 - a_3 \leq 2 \\
&1 \leq a_1 - a_3 - a_5 \leq 2 \\
&1 \leq a_1 - a_2 - a_4 \leq 2.
\end{align*}
Now, we look at a horizontal or vertical boundary of $X$ a bit more
carefully. If
there are true ``cuts'' on both ends then the row of $X$ next to it
has length at
least 3, a contradiction. Thus, we may assume without loss of generality that
$$
a_3 = 0.
$$
Next, we distinguish two cases.

{\it Case 1}: Assume $a_2 = 0$. Then we get $a_6 \leq 2$ and it is easy to see
that $C$ can be any curve satisfying these conditions except the
curve defined by
$(2,0,0,0,0,2)$.

{\it Case 2}: Assume $a_2 > 0$. Using the argument above for the
Eastern boundary
of $X$ we conclude
$$
a_4 = 0.
$$
We are left with two possibilities for the value of $a_5$.

{\it Case 2.1}: Assume $a_5 = a_6 - 1$. Since $a_5 < a_1 \leq a_6$ we
obtain $a_1
= a_6$. We conclude that either $a_2 = a_6 - 1$, i.e.\ $C$ is an arithmetically
Buchsbaum curve of type
$$
(k, k-1, 0, 0, k-1, k),
$$
or $a_2 = a_6 - 2$, i.e.\ $C$ is of type
$$
(k, k-2, 0, 0, k-1, k).
$$

{\it Case 2.2}: Assume $a_5 = a_6 - 2$. Since the length of the row next to the
Northern boundary of $X$ is at most 3 we get $a_2 = a_1 - 1$. It follows that
either $a_1 = a_6$, i.e.\ $C$ is isomorphic to the second kind of curve in Case
2.1, or $a_1 = a_6 - 1$, i.e.\ $C$ is of type
$$
(k, k-1, 0, 0, k-1, k+1).
$$
Summing up. we have shown that there are at most three types of minimal
tetrahedral curves satisfying Condition (\ref{eq-diam}). Since we
already know all
the curves with diameter one (by Corollary \ref{cor-HR-abm}), it
remains to show
that the diameter of the curves specified in the statement is at most two.

Let $C$ be the curve defined by
$$
(k, k-1, 0, 0, k-1, k+1).
$$
Enumerate the facets of the cell complex $X$ of $C$ by
$e_1,\ldots,e_{2 k}$ such
that $e_1$ is the facet in the Southeast corner and the facets $e_i$
and $e_{i+1}$
are neighbours. This enumeration exists and is uniquely determined due to the
shape of $X$. We will show that
\begin{equation} \label{eq-suff}
(a, b, c, d)^2 \cdot e_i \in \im \psi \quad \mbox{for all} \; i = 1,\ldots,2k.
\end{equation}
First, we consider $e_1$. Since $(a, b, d) \cdot e_1 \in \im \psi$ it
suffices to
show $c^2 e_1 \in \im \psi$. But this follows easily using
$$
c \cdot (-c e_1 + d e_2) \in \im \psi
$$
and
$$
c e_2 \in \im \psi.
$$
Next, we consider $e_2$. Since also $a e_2 \in \im \psi$ it remains
to show $(b,
d)^2 \cdot e_2 \in \im \psi$.  Using the minimal generators $-c e_1 + d e_2, b
e_1, d e_1$ we see that
$$
(bd, d^2) e_2 \in \im \psi.
$$
We also have
$$
b^2 e_2 \in \im \psi
$$
because $-b e_2 + a e_3, b e_3 \in \im \psi$.

Continuing similarly in this fashion we can show the relations (\ref{eq-suff})
that imply Condition (\ref{eq-diam}). This completes the argument.
\end{proof}

\begin{corollary} \label{cor-unobs}
Let $C$ be a curve that is isomorphic to one of the curves
\begin{align*}
&(k, k-1, 0, 0, k-1, k) \quad \mbox{where} \quad k \geq 1, \\
&(k, k-1, 0, 0, k-1, k+1) \quad \mbox{where} \quad k \geq 1, \quad \mbox{or} \\
&(k, k-2, 0, 0, k-1, k) \quad \mbox{where} \quad k \geq 2.
\end{align*}
Then $C$ is unobstructed and  the corresponding component of the
Hilbert scheme is
generically smooth and has dimension $4 \cdot \deg C$.
\end{corollary}

\begin{proof}
This follows by combining Proposition \ref{lem-dolcet}, Lemma
\ref{lem-diam-two},
and Theorem \ref{thm-res}.
\end{proof}

\begin{remark}
It was proved in \cite{BM7} that if $C$ is a curve in $\mathbb P^3$ (not
necessarily tetrahedral) with {\em natural cohomology} and if the diameter of
$M(C)$ is two, then $C$ is unobstructed and its Hilbert scheme has the expected
dimension.  ``Natural cohomology" means that for any $t$, no two of
$h^0(\mathcal
I_C (t))$, $h^1(\mathcal I_C(t))$ and $h^2(\mathcal I_C(t))$ are
non-zero.  Note that minimal tetrahedral curves whose Hartshorne-Rao modules
have diameter two have natural cohomology.

Note also
that our Proposition \ref{lem-dolcet}, while similar, is  independent
of that result.  For example, a non-minimal tetrahedral curve with linear
resolution (e.g.\ $(6,5,1,0,4,6)$) does not have natural cohomology, but is
unobstructed by Proposition \ref{lem-dolcet}.  On the other hand, a minimal
arithmetically Buchsbaum curve $C$ with $M(C)$ of diameter two does not have a
linear resolution (since the dual module would then be generated in
more than one
degree), but it does have natural cohomology (cf.\ \cite{BM2}, for instance
Corollary 2.5 with $t=2$ and $h=0$) and so is unobstructed.

We also remark that Mir\'o-Roig has given sufficient conditions on
the numerical
character of an {\em irreducible} arithmetically Buchsbaum curve $C$ of maximal
rank with $M(C)$ of diameter one or two, for $C$ to be unobstructed (cf.\
\cite{MMR}).
\end{remark}

The condition on the diameter of the minimal tetrahedral curves in Proposition
\ref{lem-dolcet} is sufficient for unobstructedness, but not necessary. This is
shown by the curves \linebreak $(a_1,0,0,0,0,a_6)$ because their
diameter is $a_1 +
a_6 - 1$, but all these curves are unobstructed. This follows from
the following
result.

\begin{proposition} \label{prop-n-mod}
The global sections of the normal sheaf of the curve $C$ defined by
$(a_1,0,0,0,0,a_6)$ are given by $$ H^0_*(\cN_C) \cong (R/(a,
b))^{a_1 (a_1 + 1)}
(1) \oplus (R/(c, d))^{a_6 (a_6 + 1)} (1).
$$
\end{proposition}

Before turning to the proof of this statement we record its announced
consequence.

\begin{corollary} \label{cor-skew}
Let $C \subset \PP^3$ be the curve defined by $(a_1,0,0,0,0,a_6)$.
Then its normal
sheaf satisfies
$$
H^1(\cN_C(-2)) = 0.
$$
Therefore, the corresponding component of the Hilbert scheme is generically
smooth of dimension $4 \cdot \deg C$.
\end{corollary}

\begin{proof}
The claim easily follows by Proposition \ref{prop-n-mod} and duality.
\end{proof}

The proof of Proposition \ref{prop-n-mod} is based on the computation of the
normal module of  infinitesimal neighbourhoods of a line. Recall that
the normal
module of the curve $C$ is
$$
N_C = \Hom(I_C, R/I_C).
$$

\begin{lemma} \label{lem-line} For any positive integer $k$ there is a graded
isomorphism
$$
\Hom((a, b)^k, R/(a, b)^k) \cong (R/(a, b))^{k (k+1)} (1).
$$
\end{lemma}

\begin{proof}
Put $I = (a, b)^k$ and $A = R/I$. The minimal free resolution of $I$ is
well-known
$$
0 \to R^k (-k-1) \to R^{k+1} (-k) \to I \to 0.
$$
Dualizing it with respect to $R$ provides the minimal free resolution of the
canonical module $K_A$ of $A$
$$
0 \to R \to R^{k+1} (k) \to R^k (k+1) \to K_A (4) \to 0.
$$
Dualizing with respect to $A$ we get the exact sequence
$$
0 \to \Hom(I, A) \to A^{k+1} (k) \stackrel{\psi}{\longrightarrow} A^k (k+1) \to
\Ext^1(I, A) \to 0.
$$
Thus, we see in particular that $\Ext^1 (I, A) \cong K_A \otimes A (4)$ which
allows to compute its Hilbert function. Using the last exact sequence above we
obtain
$$
\rank_K [\Hom(I, A)]_j = \left \{ \begin{array}{ll}
0 & \mif j \leq -2\\
k (k+1) (j+2) & \mif j \geq -1
\end{array}  \right.
$$
Now, let $\{e_1,\ldots,e_{k+1}\}$ be the canonical basis of the free $A$-module
$A^{k+1}$. Using $(a, b) \cdot (a, b)^{k-1} = I$ we see that
$$
(a, b)^{k-1} e_i \subset [\ker \psi]_{-1} \quad \mbox{for all} \; i =
1,\ldots,k+1.
$$
Thus, the $R$-module $M$ that is generated by
$$
G := \{a^j b^{k-1-j} e_i \s j=0,\ldots,k-1, \; i = 1,\ldots,k+1\}
$$
is a submodule of $\Hom(I, A)$. It is not too difficult to see that the minimal
generators of $M$ are annihilated by the ideal $(a, b)$ and, thus,
$G$ is a basis
of the free $R/(a, b)$-module $M$. Comparing Hilbert functions we
conclude that $M
= \Hom(I, A)$ completing the proof.
\end{proof}

We are ready for the proof of Proposition \ref{prop-n-mod}.

\begin{proof}[Proof of Proposition \ref{prop-n-mod}]
Let $C_1, C_2$ be the curves defined by $(a, b)^{a_1}$ and $(c, d)^{a_6}$,
respectively. Using the exact sequence
$$
0 \to I_C \to I_{C_1} \oplus I_{C_2} \to I_{C_1} + I_{C_2} \to 0
$$
we get
$$
H^0_*(\cO_C) \cong H^0_*(\cO_{C_1}) \oplus H^0_*(\cO_{C_2}).
$$
Since $H^0_*(\cN_C) \cong \Hom(I_C, H^0_*(\cO_C))$ it is not too
difficult to see
that the claim follows by Lemma \ref{lem-line}.
\end{proof}


\section{Remarks and problems} \label{sec-problems} 

In this section we collect some observations that do not quite fit into earlier
sections of this paper, as well as some natural questions that arise from this
work.  The list of questions is rather long, highlighting the 
richness of this line
of inquiry. The authors plan to continue investigating these questions.

\begin{remark}\label{how many}

It is natural to ask, among all tetrahedral curves, ``how many"
are minimal in their even liaison class.  That is, can we describe 
the density of
minimal tetrahedral curves among all tetrahedral curves?  The first answer is
``almost none," since of course every such even liaison class has 
infinitely many
tetrahedral curves, but essentially one minimal one (but see question 
\ref{how many
min}.\ below).  But paradoxically, we can take a different point of 
view that shows
that there are more than one might think.
\end{remark} 

We want to investigate how many minimal tetrahedral curves are in a finite set 
of tetrahedral curves. To this end it seems reasonable to fix the maximum
entry of the 6-tuple. After a change of coordinates we may assume that this
entry is the last one. Then we get:

\begin{lemma} \label{cor-how-many} 
The number of minimal tetrahedral curves $(a_1,\ldots,a_6)$ such that $a_6 = 
max \{a_1,\ldots,a_6\}$ is 
$$
N(a_6) := \sum_{a_1 = 0}^{a_6} \sum_{a_2 = 0}^{a_1 - 1} \sum_{a_5 = 0}^{a_1 -
1} \min 
\{a_1 - a_5, a_6 - a_2\} \cdot \min \{ a_1 - a_2, a_6 - a_5\}. 
$$
\end{lemma}  

\begin{proof} 
This is a consequence of Corollary \ref{s-min}. If $(a_1,\ldots,a_6)$ is 
minimal then we have $0 \leq a_1 \leq a_6$ and $0 \leq a_2, a_5 < a_1$.
Moreover, having chosen $a_1, a_2, a_5$, the only condition for $a_3$ and $a_4$
is $0 \leq a_3 < \min \{a_1 - a_5, a_6 - a_2\}$ and $0 \leq a_4 < \min \{ a_1
- a_2, a_6 - a_5\}$, respectively. The claim follows. 
\end{proof} 

The number of tetrahedral curves with fixed $a_6 = \max \{a_1,\ldots,a_6\}$ is 
$(a_6 + 1)^5$. The lemma shows that the number of minimal curves among them is
also of order $a_6^5$. To see this it suffices to find a lower estimate of the
number $N(a_6)$.  Since both minima in the formula in Lemma \ref{cor-how-many} 
are at least $a_1 - \min \{a_2, a_5\}$ we get 
\begin{eqnarray*}
N(a_6) & \geq & \sum_{a_1 = 0}^{a_6} \sum_{a_2 = 0}^{a_1 - 1}  \sum_{a_5 =
0}^{a_1 - 1} \left [ a_1 - \min \{a_2, a_5\} \right ]^2 \\ & = & \sum_{a_1 =
0}^{a_6} \sum_{a_2 = 0}^{a_1 - 1} \left [\sum_{a_5 = 0}^{a_2 - 1} (a_1 -
a_5)^2 + \sum_{a_5 = a_2}^{a_1 - 1} (a_1 - a_2)^2 \right ] \\  & = & \sum_{a_1
= 0}^{a_6} \sum_{a_2 = 0}^{a_1 - 1} \left [\sum_{k = a_1-a_2+1}^{a_1} k^2 \; +
\;  (a_1 - a_2)^3 \right ] 
\end{eqnarray*}
Now it is easy to see that this lower estimate of $N(a_6)$ is a polynomial in 
$a_6$ of order 5. In other words, we get the somewhat surprising result:  The
probability to find a minimal curve among the set of  tetrahedral curves
$(a_1,\ldots,a_6)$ such that $a_6 = \max \{a_1,\ldots,a_6\}$ is positive,
even as $a_6 \rightarrow \infty$. 

\begin{remark}\label{other cis}
The results in this paper can be extended as follows.  Let $(F_1,F_2,F_3,F_4)$
be a regular sequence.  Define curves by $I_{C_1} = (F_1,F_2), 
I_{C_2} = (F_1,F_3),
\dots, I_{C_6} = (F_3,F_4)$, and let $I_C = I_{C_1} \cap \dots \cap I_{C_6}$.
Taking $(F_1,F_2,F_3,F_4) = (a,b,c,d)$ gives the context of this 
paper. As before,
we denote by $(a_1,a_2,a_3,a_4,a_5,a_6)$ the ideal $I_C$. Note that 
no two of the
curves
$C_i$ can have a common component.  

Many of our results extend to this situation.  The key observation is 
that we still
have
\[
F_i \cdot (F_i,F_j)^{n-1} + (F_j^n) = (F_i,F_j)^n,
\]
as we noted at the beginning of the proof of Proposition \ref{reduce}.  Then
Proposition \ref{reduce} and Corollary \ref{s-min} (and indeed virtually all of
section 3) are still true, as stated, in this more general context.  Turning to
Section 4, Theorem \ref{thm-res} is clearly not true as stated.  If 
the $F_i$ all
have the same degree, it is fairly easy to modify the statement, and 
in fact the
linearity of the resolution is preserved.  If the $F_i$ have 
different degrees, it
can still be modified, but the resolution is no longer linear.  The 
other results
are less obviously generalizable to this context.
\end{remark}

We end with some open questions.  

\begin{question}
\begin{enumerate}
\item \label{how many min}
Fix an even liaison class containing tetrahedral curves.  We know 
that among the
minimal curves there is at least one that is tetrahedral.  
%
  When is a minimal tetrahedral
curve the {\em unique} minimal curve in the even liaison class (see Remark
\ref{rem-mini}, (iii))?  Note that sometimes two curves that are projectively
equivalent are linked, and sometimes they are not!
\\

\item Is it possible to characterize the tetrahedral curves with linear 
resolutions?  Recall
that this set of curves contains more than just the minimal tetrahedral curves
(Remark~\ref{rem-mini} (ii)).\\

\item We have seen that several of the minimal tetrahedral curves are 
unobstructed.
Computer experiments show that there are more of them than the ones described
above.  Are all minimal tetrahedral curves unobstructed?  Are all 
tetrahedral curves with linear resolution
unobstructed? Are all tetrahedral
curves unobstructed?
\\

\item How ``dense" are our minimal curves that give rise to nice components of 
the Hilbert scheme $H_{d, g}$, i.e.\  components that are generically smooth of
dimension $4 d$? In other words: How big are the ``gaps" in the sets of pairs
$(d,g)$  such we cannot find such a minimal curve with that degree and genus? 
\\

\item Can \ the tetrahedral curves in $\proj{3}$ that are 
arithmetically
Cohen-Macaulay be identified by explicitly giving the 6-tuples (as Schwartau 
does  for 4-tuples)?\\

\item Is there a combinatorial description of the minimal free resolution of 
the deficiency module of a tetrahedral curve?  Can we at least express the
dimensions of the components using the entries of the 6-tuple? \\

\item Can the same kind of program be carried out  in higher
projective space? Now the new question of
local Cohen-Macaulayness arises.  This includes codimension two cases 
and higher
codimension cases.  The former can still rely on complete intersection liaison
techniques, but the latter will require Gorenstein liaison techniques.
\\

\item On a computer algebra program such as {\tt macaulay} 
\cite{macaulay} one can
make experiments, and one notices a great deal of structure in the 
Betti diagram.
In particular, non-minimal tetrahedral curves have linear strands in the
resolution.   How do these arise? Certainly they come from the Betti 
diagram of the
minimal tetrahedral curve using basic double links, but this should 
be explained in
a clearer way.  Can one relate the ``degree of non-minimality'' to 
the number of
linear strands?  The result of \cite{MN1} Corollary 4.5 should be useful.\\

\item Is it possible to use similar techniques to study more general monomial 
ideals?

\end{enumerate}
\end{question}


\section*{Appendix: The algorithm to find $S$-minimal curves}

Below we record a crude MAPLE implementation of Algorithm \ref{alg}. 
It is a MAPLE
work sheet that can be downloaded at 
\medskip
{\small 
\begin{center}
{\tt http://www.ms.uky.edu/$\sim$uwenagel/}.
\end{center} }

\bigskip

\begin{verbatim}
>  # algorithm to compute the S-minimal curve of a tetrahedral curve
   # Input: weight vector of the given curve

>  smin := proc(aa,bb,c,d,e,f)
   local a, b, m, i, j, k, s, t, w, W, x, y, z, output1, output2;

   a := [aa,bb,c,d,e,f];
   b:= a; s := 0;

   if (s > -1) then do

   # 1. computation of the maximal weight m of a line
   m := a[1]; j := 1;
   for i from 2 to 6 do
                        if (a[i]  > m) then m := a[i];
                                            j := i;
                        end if;
                     od;

   # 2. avoid unnecessary computations
   if (m = 0) then
        output1 := matrix(1,6,[[` `,`Minimal curve to`, b, `is`, a,` `]]):
        output2:= matrix(1,3,[[`It is obtained after`, s,` reduction(s).`]]):
        print(output1); print(output2);
        return(a);
   end if;

   # 3. computation of facet of maximal weight W
   x := a[1] + a[2] + a[3];
   y := a[1] + a[4] + a[5];
   z := a[2] + a[4] + a[6];
   t := a[3] + a[5] + a[6];

   w := [x, y, z, t];

   W:= w[1]; k := 1;

   for i from 2 to 4 do
                        if (w[i]  > W) then W := w[i];
                                            k := i;
                        end if;
                     od;


   # 4. test if curve is S-minimal
   if (a[j] + a[7-j] > W) then
      output1 := matrix(1,6,[[` `,`Minimal curve to`, b, `is`, a,` `]]):
      output2:= matrix(1,3,[[`It is obtained after`, s,` reduction(s).`]]):
      print(output1); print(output2);
      return(a);
   end if;

   # 5. reduction of the non-minimal curve
   s := s+1;

   if (k = 1) then a[1] := max(0, a[1] - 1);
                   a[2] := max(0, a[2] - 1);
                   a[3] := max(0, a[3] - 1);

   end if;

   if (k = 2) then a[1] := max(0, a[1] - 1);
                   a[4] := max(0, a[4] - 1);
                   a[5] := max(0, a[5] - 1);

   end if;

   if (k = 3) then a[2] := max(0, a[2] - 1);
                   a[4] := max(0, a[4] - 1);
                   a[6] := max(0, a[6] - 1);

   end if;

   if (k = 4) then a[3] := max(0, a[3] - 1);
                   a[5] := max(0, a[5] - 1);
                   a[6] := max(0, a[6] - 1);

   end if;

   end do;

   end if;

   end:

\end{verbatim}
\bigskip

In order to use the procedure above in a different work sheet one 
could write it into a text file named, for example, proc.txt. Then it 
can be used as follows.

\begin{verbatim}

>  read `C:\\suitable path\\proc.txt`;
>  smin(5,1,3,2,2,5):

         [  , Minimal curve to , [5, 1, 3, 2, 2, 5] , is ,

         [5, 1, 2, 2, 1, 4] ,  ]


             [It is obtained after    1     reduction(s).]

>  smin(6,0,8,1,0,4):

         [  , Minimal curve to , [6, 0, 8, 1, 0, 4] , is ,

         [0, 0, 0, 0, 0, 0] ,  ]


             [It is obtained after    10     reduction(s).]
\end{verbatim}

\end{document}